\documentclass[reqno,11pt]{amsart}
\usepackage{amssymb}
\usepackage{amsmath}
\usepackage{amsthm}
\usepackage{pgf}
\usepackage{color}
\usepackage{mathrsfs}
\usepackage{graphicx}
\usepackage{hyperref}
\usepackage{graphicx}
\usepackage{tikz}
\usepackage{pgf}
\usepackage{enumerate}
\usepackage{ulem}
\usepackage{mathtools}

\usepackage{bm}

\usepackage{enumitem}

\usepackage{comment}

\newtheorem{theorem}{Theorem}[section]
\newtheorem{lemma}[theorem]{Lemma}

\newtheorem{remark}[theorem]{Remark}

\renewcommand{\lvert}{\|}
\renewcommand{\rvert}{\|}

\newcommand{\grad}{\bm{\nabla}}

\DeclareMathOperator{\spn}{span}

\newcommand{\mathsout}[1]
{\bgroup\mathchoice
	{\sbox0{$\displaystyle{#1}$}%
		\usebox0\hspace{-\wd0}%
		\rule[0.5\ht0-0.5\dp0-.5pt]{\wd0}{0.1pt}}%
	{\sbox0{$\textstyle{#1}$}%
		\usebox0\hspace{-\wd0}%
		\rule[0.5\ht0-0.5\dp0-.5pt]{\wd0}{0.1pt}}%
	{\sbox0{$\scriptstyle{#1}$}%
		\usebox0\hspace{-\wd0}%
		\rule[0.5\ht0-0.5\dp0-.5pt]{\wd0}{0.1pt}}%
	{\sbox0{$\scriptscriptstyle{#1}$}%
		\usebox0\hspace{-\wd0}%
		\rule[0.05\ht0-0.5\dp0-.5pt]{\wd0}{0.1pt}}%
	\egroup}

\renewcommand{\d}{{\mathrm d}}

\renewcommand{\div}{\operatorname{div}}

\def\to{\rightarrow}

\newtheorem{definition}[theorem]{Definition}

\numberwithin{equation}{section}

\makeatletter
\newcommand*\bigcdot{\mathpalette\bigcdot@{.5}}
\newcommand*\bigcdot@[2]{\mathbin{\vcenter{\hbox{\scalebox{#2}{$\m@th#1\bullet$}}}}}
\makeatother

\begin{document}

\title[Analysis of a Navier-Stokes Phase-Field Crystal system]
{Analysis of a Navier-Stokes \\ Phase-Field Crystal system}
\author[C. Cavaterra, M. Grasselli, M.A. Mehmood, R. Voso]
{\textsc{Cecilia Cavaterra}$^\dagger$,
\textsc{Maurizio Grasselli}$^\ast $\\
\textsc{Muhammed Ali Mehmood}$^\ddagger$
\& \textsc{Riccardo Voso}$^\sharp$}

\address{$^\dagger$ Dipartimento di Matematica ``F. Enriques''\\
Universit\`{a} degli Studi di Milano\\
20133 Milano, Italy \& IMATI-CNR 27100 Pavia, Italy}
\email{cecilia.cavaterra@unimi.it}
\address{$^\ast$Dipartimento di Matematica\\
Politecnico di Milano\\
20133 Milano, Italy}
\email{maurizio.grasselli@polimi.it}
\address{$^\ddagger$Department of Mathematics\\
Imperial College London\\
London, SW7 2AZ, UK}
\email{muhammed.mehmood21@imperial.ac.uk}
\address{$^\sharp$ Vienna School of Mathematics\\
University of Vienna\\
1090 Wien, Austria}
\email{riccardo.voso@univie.ac.at}



\allowdisplaybreaks

\begin{abstract}
We consider an evolution system modeling a flow of
colloidal particles which are suspended in an incompressible fluid
and accounts for colloidal crystallization. The system consists of
the Navier-Stokes equations for the volume averaged velocity
coupled with the so-called Phase-Field Crystal equation for the density deviation.
Considering this system in a periodic domain and assuming that the viscosity as well as the mobility depend on the density deviation, we first prove
the existence of a weak solution in dimension three. Then, in dimension two,
we establish the existence of a (unique) strong solution.

\end{abstract}

\maketitle

\thispagestyle{empty}

\section{Introduction}
In this paper, we consider a hydrodynamic phase-field crystal (PFC) model introduced in \cite{nspfc} (see also \cite{Praetorius,Van} and their references) to describe the colloidal suspension in a fluid, providing, in particular, a quantitative approach down to the characteristic length scale of the colloidal particles and accounting for colloidal crystallization. Physical processes that lead to permanently evolving pattern, such as Rayleigh-B\'enard convection, nucleation, and crystal growth, are tackled using the Dynamical Density Functional Theory approach. More precisely, the phase-field variable (or order parameter) is the density deviation $\phi$ whose associated free energy is the so-called Swift-Hohenberg (SH) free-energy functional (see \cite{Swift}, see also \cite{SH}). In a dimensionless form, this is given by
\begin{equation*}
\mathcal{F}_{\text{sh}}(\phi)=\int_Q  \left(   \frac{1}{2} \vert \Delta \phi \vert^2 - q_0^2 \vert \nabla \phi \vert^2 + \frac{1}{4}\phi^4 + \frac{r + q_0^4}{2} \phi^2    \right)\, \d x,
\end{equation*}
where $Q = (0,L)^d$, $d=2,3$, $L>0$, is a given periodic domain.
Here, $r\in \mathbb{R}$ is a phenomenological parameter and the constant $q_0$ is related to the lattice spacing and is usually set equal to the unity. The evolution of $\phi$ is thus governed by the so-called PFC equation, that is, a conserved SH equation (see, for instance, \cite{Backofen,Emmerich})
\begin{align*}
&\partial _t \phi + \bm{u} \cdot \nabla \phi = -\div (m(\phi) \nabla \psi ),\\
&\psi=\frac{\delta \mathcal{F}_{\text{sh}}[\phi]}{\delta \phi},
\end{align*}
in $Q\times (0,T)$, $T>0$. Here, $m(\cdot)$ is the non-costant mobility coefficient, while $\bm{u}$ stands for the
(volume) averaged fluid velocity whose evolution is governed by
\begin{align*}
&\partial _{t} \bm{u} + \bm{u} \cdot \nabla \bm{u} = \div (\eta (\phi) \bm{D} \bm{u}) + \nabla P - M \phi \nabla \psi,\\
&\div \bm{u} = 0,
\end{align*}
where $\eta(\cdot)$ is the variable viscosity, $\bm{D} \bm{u}$ is the usual shear strain rate tensor, $M>0$ is a given constant, and other constants have been set equal to the unity.

Summing up, we have the following Navier-Stokes Phase-Field Crystal (NSPFC) system (see \cite[(28)]{nspfc})
\begin{align}
&\partial _{t} \bm{u} + \bm{u} \cdot \nabla \bm{u} = \div (\eta (\phi) \bm{D} \bm{u}) + \nabla P - M \phi \nabla \psi,\label{eq:ns}\\
&\div\bm{u} = 0, \label{eq:div}\\
&\partial _t \phi + \bm{u} \cdot \nabla \phi = \nabla \cdot (m(\phi) \nabla \psi),\label{eq:ch}\\
&\psi = \Delta^2 \phi +2\Delta \phi + f(\phi), \label{eq:F}
\end{align}
subject to the initial and the periodic boundary conditions
\begin{align}
&\bm{u}(\cdot,0)=\bm{u}_0, \quad \phi(\cdot,0)=\phi_0, \quad \text{ in } Q, \label{eq:ci}\\
&\bm{u}(x,t)=\bm{u}(x +L\bm{e}_i,t) , \quad
\phi(x,t)=\phi(x+L\bm{e}_i,t), \quad x\in \mathbb{R}^d, \quad i=1,...,d ,\label{eq:periodic}
\end{align}
where $(\bm{e}_i)_{i=1}^{d}$ is the canonical basis of $\mathbb{R}^d$ and $\bm{u}_0$, $\phi_0$ are given initial data. Here $$ f(\phi) = \phi^3 + (r + 1) \phi, $$ having set $q_0=1$. Moreover, we assume that $\eta$ and $m$ are
smooth enough and bounded from above and below by suitable positive constants (see next section for details).

The main goal of this work is to analyze the existence and, in dimension two, the uniqueness of
solutions to the above problem \eqref{eq:ns}-\eqref{eq:periodic}.
More precisely, we first establish the existence of a (global) weak solution which is unique in dimension two even in the case of non constant viscosity (cf. \cite{AGiorg} and its references for the Navier-Stokes-Cahn-Hilliard system). Also, in dimension two, we prove the existence of a (unique) strong solution which continuously depends on the initial data. Besides, we prove that any weak solution becomes strong in finite time.

There are several interesting open issues that might be explored in a future paper. For instance, the existence of a local strong solution in dimension three, the asymptotic convergence of a solution towards a single stationary state, and the existence of (global and, possibly, exponential) attractors within the theory of infinite-dimensional dissipative dynamical systems. Also, the present analysis could be extended to a binary PFC model (see \cite{Yang-He_2022}) or to a model where the PFC equation is replaced with a conserved SH equation (see \cite{Lee}) or with the so-called modified PFC equation (see, for instance, \cite{GW1,GW2} and references therein). However, in the latter case, the analysis seems rather challenging.

The plan of the paper goes as follows. The next section is
divided in four subsections. The first two are
devoted to introduce the functional framework and to recall some basic tools. Then the notions of weak and strong solutions are defined and, in the last subsection, the main results are stated. Section 3 contains the proof of the existence of a (global) weak solution. The existence of a strong solution and its continuous dependence on the initial data as well as the regularization of any weak solution are established in Section~4.

\section{Preliminaries}
\subsection{Notation and functional spaces}
Let $B$ be a (real) Banach space. We denote by $\lVert \cdot \rVert_B$ and by $\langle \cdot , \cdot \rangle_{B',B}$ its norm and the duality product between $B$ and its topological dual $B'$, respectively. We indicate by $\mathcal{L}(X,Y)$ the space of linear bounded operators from a (real) Banach space $X$ into another one $Y$ and we set $\mathcal{L}(X):=\mathcal{L}(X,X)$. For all functions $w \in L^1(Q)$, $\langle w \rangle := \vert Q \vert^{-1}_d \int_Q w(x) \,\d x$ denotes the mean value of $w$. Here $\vert \cdot\vert_d$ stands for the Lebesgue measure. Without loss of generality, in the definition of $Q$ we can take $L=1$ so that $\vert Q \vert_d  =1$. Also, in the sequel, we shall denote by $c$ a generic positive constant which may depend on $T$ and, possibly, on other given quantities.

Let $C_p^{\infty}(Q)$ be the space of functions in $C^{\infty}(Q)$ which are periodic over $Q$. For any $s \in \mathbb{N}$, we indicate by $\Phi_s$ the Sobolev spaces of functions in  $H^s(Q)$ which are periodic over $Q$, namely the completion of $C_p^{\infty}(Q)$ with respect to the $H^s(Q)$-norm $$\Phi_s=\overline{\{ w \in C_p^{\infty}(Q)\} }^{H^s(Q)}.  $$ These spaces are Hilbert spaces with respect to the scalar product
$$(w,v)_{H^{s}(Q)} = \sum_{\vert k \vert \leqslant s} (D^k w, D^k v)_{L^2(Q)}.$$
Therefore, the induced norm is   $\|w\|_{H^{s}(Q)}= \sqrt{(w,w)_{H^{s}(Q)}}$.

Let $\mathbb{B} = B^d$ endowed with the product structure indicate the space of vector-valued functions $\bm{w} := (w_1,...,w_d)$, namely $\bm{w} \in \mathbb{B}$ if and only if $w_i \in B$ for every $i =1,...,d$. We introduce the following Hilbert spaces $$ \mathbb{H}=\overline{\{ \bm{w} \in \mathbb{C}_p^{\infty}(Q) \; : \; \div \bm w =0 ,\; \langle \bm{w} \rangle = 0\} }^{\mathbb{L}^2(Q)}, $$ $$ \mathbb{V}=\overline{\{ \bm{w} \in \mathbb{C}_p^{\infty}(Q) \; : \; \div \bm w =0 ,\; \langle \bm{w} \rangle = 0\} }^{\mathbb{H}^1(Q)}. $$ Scalar product and induced norm in $\mathbb{H}$ are then defined in the canonical sense, namely $(\bm{w},\bm{v})_{\mathbb{H}} = \sum_{ i= 1}^d ( w_i,v_i)_{L^2(Q)} $ and $\lvert \bm{w} \rVert_{\mathbb{H}} = \sqrt{(\bm{w},\bm{w})_{\mathbb{H}}}$, respectively, and scalar product and the induced norm in $\mathbb{V}$ are given by $(\bm{w},\bm{v})_{\mathbb{V}} = \sum_{ i= 1}^d ( \nabla w_i, \nabla v_i)_{\mathbb{H}}$ and $\lvert \bm{w} \rVert_{\mathbb{V}} = \sqrt{(\bm{w},\bm{w})_{\mathbb{V}}}$, respectively. As the Poincar\'e inequality $$ \lvert \bm{w} \rvert_{\mathbb{H}} \leq c \lvert \nabla \bm{w} \rvert_{\mathbb{H}}, \quad \forall \bm{w} \in \mathbb{V}  $$ holds, the $\mathbb{V}$-norm is equivalent to the canonical $\mathbb{H}^1(Q)$-norm on $\mathbb{V}$ and, thanks to the Korn equality (see, for instance, \cite[Lemma 1.2, p.4]{BoyerArticle}), i.e. $\lvert \nabla \bm{w}_i \rvert_{\mathbb{H}} = \sqrt{2} \lvert \bm{D}\bm{w}_i \rvert_{\mathbb{H}} $ for all $\bm{w} \in \mathbb{V}$, another equivalent norm on $\mathbb{V}$ is given by $\lvert \bm{D} \bm{w} \rvert_{\mathbb{H}}$. In the following, we consider two distinct Hilbert triplets $\mathbb{V} \hookrightarrow \mathbb{H} \cong \mathbb{H}' \hookrightarrow \mathbb{V}'$ and $\Phi_4 \hookrightarrow \Phi_2 \cong (\Phi_2)' \hookrightarrow (\Phi_4)'$, with the usual identifications $\mathbb{H} \cong \mathbb{H}'$ and $\Phi_2 \cong (\Phi_2)'$. We recall that, the Rellich-Kondrachov theorem ensures that $\mathbb{V}$ and $\Phi_4$ are compactly embedded in $\mathbb{H}$ and $\Phi_2$, respectively.

Eventually, we recall the following version of Poincar\'e inequality which will be used throughout the paper (see, e.g., \cite[Lemma 1.5, p.5]{BoyerArticle})
\begin{lemma}
\label{poincphi}
For any $w \in \Phi_1$ we have $$\lVert w - \langle w \rangle \rVert_{\Phi_1} \leqslant c \lVert \grad w \rVert_{L^2(Q)},$$ and for any $ w \in \Phi_{s+2} $ we have $$\lVert w - \langle w \rangle \rVert_{\Phi_{s+2}} \leqslant c  \lVert \Delta w \rVert_{\Phi_s}.$$
\end{lemma}
In particular, Lemma \ref{poincphi} entails that $\lVert \grad w \rVert_{\Phi_1} \leqslant c \lVert \Delta w \rVert_{L^2(Q)}$ for any $w \in \Phi_2$.

\subsection{Inertia Term}
We recall that the trilinear form
\begin{equation*}
b_0(\bm{u},\bm{v},\bm{w})= \int_Q \bigl((\bm{u} \cdot \nabla\bm{v})\cdot\bm{w}\bigr), \qquad \forall \bm{u},\bm{v},\bm{w} \in \mathbb{V}
\end{equation*}
is continuous on $\mathbb{V}\times\mathbb{V}\times\mathbb{V}$ and satisfies the identities
\begin{align}
	b_0(\bm{u},\bm{v},\bm{w}) + b_0(\bm{u},\bm{w},\bm{v}) &=0, \qquad \forall \bm{u},\bm{v},\bm{w} \in \mathbb{V},\\
	b_0(\bm{u},\bm{v},\bm{v}) &=0, \qquad \forall \bm{u},\bm{v}\in \mathbb{V}.
\end{align}
By using the Ladyzhenskaya's inequality it can be shown that, for all $\bm{u},\bm{v},\bm{w} \in \mathbb{V}$,
\begin{equation}
\label{estimateb}
	\lvert b_0 (\bm{u},\bm{v},\bm{w}) \rvert \leq c \lVert \bm{u} \rVert_{\mathbb{H}}^{1 - d/4}  \lVert \bm{u} \rVert_{\mathbb{V}}^{d/4} \lVert \bm{v} \rVert_{\mathbb{H}}^{1- d/4}  \lVert \bm{v} \rVert_{\mathbb{V}}^{d/4}  \lVert \bm{w} \rVert_{\mathbb{V}}, \qquad d=2,3.
\end{equation}
Then we can introduce the {\it inertia term} $B_0$, namely the continuous bilinear operator from $\mathbb{V}\times \mathbb{V}$ to $\mathbb{V}'$ defined as
$$
	\langle B_0(\bm{u},\bm{v}) , \bm{w} \rangle = b_0(\bm{u},\bm{v},\bm{w}).
$$
Thus, from \eqref{estimateb}, we readily get
\begin{equation}
\label{estimateB}
	\lVert B_0(\bm{u},\bm{u}) \rVert_{\mathbb{V'}} \leq c \lVert \bm{u} \rVert_{\mathbb{H}}^{2 - d/2} \lVert \bm{u} \rVert_{\mathbb{V}}^{d/2}, \qquad d=2,3, \ \forall \bm{u} \in \mathbb{V}.
\end{equation}

\subsection{Definitions of Solution}
The aim of this section is to give a rigorous definition of weak and strong solution to problem \eqref{eq:ns}-\eqref{eq:periodic}. The first one goes as follows

\begin{definition}[Weak solutions]
\label{definitionWeak}
Let $\bm{u}_0\in \mathbb{H}$, $\phi_0\in\Phi_2$ be given. We say that $(\bm{u},\phi)$ is a weak solution to problem \eqref{eq:ns}-\eqref{eq:periodic} in $(0,T)$, if $\bm{u}$ and $\phi$ are such that
\begin{equation} \label{deboleteorema5}
    \bm{u}\in L^{\infty}(0,T;\mathbb{H})\cap L^2(0,T;\mathbb{V})\cap C_w([0,T];\mathbb{H}),
\end{equation}
\begin{equation} \label{deboleteorem6}
     \phi\in C^{0}([0,T];\Phi_2)\cap L^2(0,T;\Phi_5)\cap H^1(0,T;\Phi_1^\prime),
\end{equation}
and the following identities hold almost everywhere in $(0,T)$,
\begin{equation} \label{deboleteorema1}
    	( \bm{u}' , \bm{v} )  + ( \eta(\phi) \bm{D}\bm{u},\nabla \bm{v} )+ b_0(\bm{u},\bm{u},\bm{v}) + M ( \phi \nabla \psi,\bm{v} ) =0  ,
\end{equation}
\begin{equation} \label{deboleteorema2}
    ( \phi',\rho ) + ( \bm{u} \nabla \phi,\rho ) + ( m(\phi)\nabla \psi, \nabla \rho ) =0,
\end{equation}
for all $\bm{v}\in \mathbb{V}$ and for all $\rho\in \Phi_1$,
where $\psi \in L^2(0,T;\Phi_1)$ satisfies, almost everywhere in
$Q \times (0,T) $,
\begin{equation}
\label{deboleteorema3}
	\psi = \Delta^2\phi + 2\Delta \phi +f(\phi),
\end{equation}
along with the initial conditions
\begin{equation}
	\label{deboleteorema4}
	\bm{u}(0) = \bm{u}_0, \qquad \phi(0) = \phi_0.	
\end{equation}
\end{definition}

The definition of strong solution is given by
\begin{definition}[Strong solutions]
\label{definitionStrong}
Let $\bm{u}_0\in \mathbb{V}$, $\phi_0\in \Phi_3$ be given. We say that $(\bm{u},\phi)$ is a strong solution to problem \eqref{eq:ns}-\eqref{eq:periodic} in $(0,T)$, if $(\bm{u},\phi)$ is a weak solution and if, in addition, $\bm{u}$ and $\phi$ have the following properties
\begin{equation} \label{forteteorema1}
    \bm{u}\in C^{0}(0,T;\mathbb{V})\cap L^2(0,T;\mathbb{H}^2(Q)),
\end{equation}
\begin{equation} \label{forteteorema2}
    \phi\in C^{0}([0,T];\Phi_3)\cap L^2(0,T;\Phi_6)\cap H^1(0,T;\Phi_0),
\end{equation}
\begin{equation} \label{forteteorema4}
\psi\in L^{2}(0,T;\Phi_2).
\end{equation}
Therefore, $\bm{u}$ and $\phi$
 also satisfy \eqref{eq:ns}-\eqref{eq:ch} almost everywhere in $Q\times (0,T)$.
\end{definition}

\subsection{Statements of the Main Results}
In this subsection we state  our two main theorems whose proofs are given in the next sections. The first establishes the existence of a weak solution to problem \eqref{eq:ns}-\eqref{eq:periodic}, whereas the second ensures that, in two dimensions, there exists a unique strong solution to problem \eqref{eq:ns}-\eqref{eq:periodic} which continuously depends on the initial data.

Let us assume
\begin{itemize}
	\item[\textbf{(A1)}]  $\eta, m  \in W^{1,\infty}(\mathbb{R};\mathbb{R}^+)$ and there exist positive constants $\eta_0, \eta_1, \eta_2, m_0, M_0, m_2$ such that
\end{itemize}
\begin{align*}
&\eta_1 \geqslant \eta (s) \geqslant \eta_0 > 0 \quad \forall s\in\mathbb{R}, \qquad\quad\; \eta_2 \geqslant \eta'(s) >0 \quad \text{for a.a.}\;s\in \mathbb{R},\\
&M_0\geqslant m (s) \geqslant m_0 > 0\quad \forall s\in\mathbb{R}, \qquad m_2\geqslant m'(s) >0 \quad \text{for a.a.}\;s\in \mathbb{R}.
\end{align*}
Then the existence of a weak solution, which is unique
in dimension two, is given by

\begin{theorem}
\label{theoremweak}
Let $\bm{u}_0\in\mathbb{H}$, $\phi_0\in \Phi_2$ be given and let $\eta$ and $m$ satisfy (A1). Then:
\begin{enumerate}[label=(\roman*)]
	\item For $d = 2,3$, there exists at least one weak solution $(\bm{u},\phi)$ to problem \eqref{eq:ns}-\eqref{eq:periodic} in the sense of Definition \ref{definitionWeak}. For $d=3$, any weak solution will additionally satisfy the following energy inequality
	\begin{align}
 \label{enineq}
 \notag
	&\frac{1}{2}\|\bm{u}(t)\|_{\mathbb{L}^{2}(Q)}^2 + \mathcal{F}_{sh}[\phi(t)] + \int_0^t \int_{Q} (2 \eta(\phi) |D( \bm{u} )|^{2} + m(\phi)|\nabla \psi |^2) \mathrm{d}x\mathrm{d}t \\
    &\le \frac{1}{2}\|\bm{u}_0\|_{\mathbb{L}^{2}(Q)}^2 + \mathcal{F}_{sh}[\phi_0], \qquad\text{ for a.e. } t\in(0,T).
	\end{align}
	\item  If $d=2$, then a weak solution is unique and satisfies the energy identity
	\begin{align}  \label{enineq2}
    \notag
	& \frac{1}{2}\|\bm{u}(t)\|_{\mathbb{L}^{2}(Q)}^2 + \mathcal{F}_{sh}[\phi(t)]  + \int_0^t \int_{Q} (2 \eta(\phi) |D( \bm{u} )|^{2} + m(\phi)|\nabla \psi |^2) \mathrm{d}x\mathrm{d}t\\
    &=  \frac{1}{2}\|\bm{u}_0\|_{\mathbb{L}^{2}(Q)}^2 + \mathcal{F}_{sh}[\phi_0], \qquad\text{ for all } t\in [0,T].
	\end{align}
\end{enumerate}
\end{theorem}
In dimension two, we can prove the existence of a (unique) strong solution. Moreover, any weak solution
instantaneously becomes a strong solution. Indeed, we have
\begin{theorem}
\label{theoremstrong}
Let $d=2$.  Suppose $\bm{u}_0\in\mathbb{V}$, $\phi_0\in \Phi_3$ and assume that $\eta$ and $m$ satisfy (A1). Then, there exists a unique global strong solution $(\bm{u},\phi)$ to problem \eqref{eq:ns}-\eqref{eq:periodic} in the sense of Definition \ref{definitionStrong}.
Moreover, this solution depends continuously on initial data, that is, if $(\bm{u}_1,\phi_1)$ and $(\bm{u}_2,\phi_2)$ are two strong solutions of the problem \eqref{eq:ns}-\eqref{eq:periodic} originated from the initial data $(\bm{u}_{0,1},\phi_{0,1})$ and $(\bm{u}_{0,2},\phi_{0,2})$ respectively, then the following inequality holds
\begin{align}
\label{dipcont}
    &\lVert\bm{u}_1-\bm{u}_2\rVert_{C^0([0,T];\mathbb{H})}^2 + \lVert\bm{u}_1-\bm{u}_2\rVert_{L^2(0,T;\mathbb{V})}^2 + \lVert\phi_1-\phi_2\rVert_{C^0([0,T];H_p^2(Q))}^2 \nonumber\\
    &+ \lVert\phi_1-\phi_2\rVert_{L^2(0,T;H_p^5(Q))}^2 \nonumber\\
    &\leqslant c \biggl( \lVert\bm{u}_{0,1}-\bm{u}_{0,2}\rVert_{\mathbb{H}}^2 + \lVert\phi_{0,1}-\phi_{0,2}\rVert_{H_p^2(Q)}^2 \biggr).
\end{align}
Moreover, any weak solution instantaneously regularizes and becomes a strong solution.
\end{theorem}

\begin{remark}
Thanks to the energy inequality \eqref{enineq}, it is possible to construct a global weak solution, i.e., a solution defined in $[0,+\infty)$ (see, e.g., \cite[Chapt.V, 1.3.6]{Boyer} for the Navier-Stokes system). This solution is unique in dimension two and it is strong in $[\tau,+\infty)$ for any fixed $\tau>0$.
\end{remark}

\begin{remark}
On account of the boundedness of $\phi$ (see \eqref{deboleteorem6}), one can take a more general function $f:\mathbb{R} \to \mathbb{R}$ of any growth, which is sufficiently smooth (e.g., $f \in C^2(\mathbb{R})$) and such that, for instance,
$$
f(y)y \geq c_1 y^4 - c_2y^2,
$$
for all $y\in\mathbb{R}$ and for some positive constants $c_1$, $c_2$.
\end{remark}

\section{Proof of Theorem \ref{theoremweak}}
This section is devoted to the proof of Theorem \ref{theoremweak}. We first establish the existence of a weak solution which satisfies the energy inequality. More precisely, in Subsection \ref{approximateweak}, a Faedo-Galerkin approximation of the problem \eqref{deboleteorema1}-\eqref{deboleteorema4} and the proof of its well-posedness are given. In Subsection \ref{estimatesweak} we provide {\it a priori} estimates on such solutions $(\bm{u}_n,\phi_n)$ which are shown, in Subsection \ref{limitweak}, to converge, up to subsequences, to a weak solution of the problem \eqref{eq:ns}-\eqref{eq:periodic}. Finally, in Subsection \ref{weakunique2D} we prove that weak solutions are unique in two dimensions as well as the validity of the energy identity.

\subsection{Approximating problem}
\label{approximateweak}
Let $(\mathbb{V}^n)_{n\geqslant1}\subset \mathbb{V}$ be a sequence of subspaces such that $\mathbb{V}^n \subseteq \mathbb{V}^{n+1}$, $\mathbb{V}^{\infty} := \bigcup_n \mathbb{V}^n$ is dense in $\mathbb{V}$ and $\dim(\mathbb{V}^n)=n<\infty$. Analogously, let $(\Phi_2^n)_{n\geqslant1} \subset \Phi_2$ be a sequence of subspaces such that $\Phi_2^n\subseteq \Phi_2^{n+1}$, $\Phi_2^{\infty} := \bigcup_n \Phi_2^n$ is dense in $\Phi_2$ and $\dim(\Phi_2^n)=n<\infty$. Let $(\bm{w}_j)_{j\geqslant1}$ be the basis made of eigenfunctions of the Stokes operator with periodic boundary conditions and let $(\rho_j)_{j\geqslant1}$ be the basis made of eigenfunctions of the operator $\Delta^2 +2\Delta$ with periodic boundary conditions (see \cite[Theorem 4.2, p.225]{Giorgini} for details). We choose $\mathbb{V}^n$ and $\Phi_2^n$ as the finite-dimensional spaces generated by $(\bm{w}_j)_{j\geqslant1}$ and $(\rho_j)_{j\geqslant1}$, respectively. Namely, $\mathbb{V}^n = \spn\{\bm{w}_1,...,\bm{w}_n\}$ and $\Phi_2^n = \spn\{\rho_1,...,\rho_n\}$. We also introduce the orthogonal projectors $P_{\mathbb{V}^n}$ and $P_{\Phi_2^n}$ from $\mathbb{H}$ in $\mathbb{V}^n$ and from $\Phi_2$ in $\Phi_2^n$, respectively. Then, set
\begin{equation*}
\bm{u}_n(\bm{x},t)=\sum_{j=1}^n a_{j,n}(t)\bm{w}_j(\bm{x}) , \quad
\phi_n(\bm{x},t)=\sum_{j=1}^n b_{j,n}(t)\rho_j(\bm{x}) , \quad \psi_n(\bm{x},t)=\sum_{j=1}^n c_{j,n}(t)\rho_j(\bm{x}),
\end{equation*}
where $a_{j,n}$, $b_{j,n}$, and $c_{j,n}$ are functions to be determined in $C^1([0,T])$, $C^1([0,T])$ and $C^0([0,T])$, respectively, such that $(\bm{u}_n,\phi_n)$ solves for all $\bm{w}\in \mathbb{V}^n$, for all $\rho\in \Phi_2^n$, and for almost every $t \in (0,T)$,
\begin{equation} \label{eq:Ia}
( \bm{u}_n' , \bm{w} )  + ( \eta(\phi_n) \bm{D}\bm{u}_n,\nabla \bm{w} )+ b_0(\bm{u}_n,\bm{u}_n,\bm{w}) + M ( \phi_n \nabla \psi_n,\bm{w} ) =0  ,
\end{equation}
\begin{equation} \label{eq:IIa}
( \phi_n',\rho ) + ( \bm{u}_n \nabla \phi_n,\rho ) + ( m(\phi_n)\nabla \psi_n, \nabla \rho ) =0,
\end{equation}
where
\begin{equation}
\label{eq:IIIa}
\psi_n = \Delta^2\phi_n + 2\Delta \phi_n +P_{\Phi_2^n}f(\phi_n),
\end{equation}
and satisfies the initial conditions
\begin{equation}
\label{eq:IIIIa}
\bm{u}_n(0) = \bm{u}_{0,n} = P_{\mathbb{V}^n}\bm{u}_0, \qquad \phi_n(0) = \phi_{0,n} = P_{\Phi_2^n} \phi_0.
\end{equation}

Existence of a unique solution $(\bm{u}_n,\phi_n) \in C^1([0,T_n); \mathbb{V}^n) \times  C^1([0,T_n),\Phi_2^n)$, $0<T_n \leq T$, of the problem \eqref{eq:Ia}-\eqref{eq:IIIIa} follows then by the Cauchy-Lipschitz-Picard theorem. We also note that $\psi_n \in C^0([0,T_n);\Phi_2^n)$. Before closing this subsection we observe that  $\phi_n \in \mathbb{C}_p^{\infty}(Q)$ so that $\phi_n \subset \Phi_s$ for any $s\in\mathbb{N}$. This turns out to be crucial in Section \ref{strong2D} when we provide estimates for the proof of the Theorem \ref{theoremstrong}.

\subsection{A priori estimates}
\label{estimatesweak}
In this subsection, we prove {\it a priori} estimates on the solutions $(\bm{u}_n,\phi_n)$ which are uniform on $n$. Here, $c$ stands for a generic positive constant which is independent of $n$, but it may depend on $T$ and on other given quantities.

$\bullet$ {\it Step 1:}
	By taking $\rho = 1$ in \eqref{eq:IIa}, we get, after integrating by parts, $\langle \phi_n' \rangle = 0$, namely
	\begin{equation}
	\label{step0}
		\langle \phi_n(t) \rangle  = \langle P_{\Phi_2^n} \phi_0 \rangle = \langle \phi_0 \rangle.
	\end{equation}
	Similarly, we test equation \eqref{eq:IIIa} with $\rho=1$. As $\rho_1 = 1$ is the eigenfunction of the operator $\Delta^2 + \Delta$ corresponding to the eigenvalue $\gamma_1 = 0$, we readily obtain $\langle \psi_n \rangle = \langle P_{\Phi_2^n} f(\phi_n) \rangle$, which implies, by exploiting the definition of $f$,
	\begin{equation}
	\label{step00}
		\vert \langle \psi_n \rangle \vert \leqslant  c \bigl( \lVert \phi_n \rVert_{L^3(Q)}^3 + \lVert \phi_n \rVert_{L^1(Q)} \bigr).
	\end{equation}

$\bullet$ {\it Step 2:}
By summing \eqref{eq:Ia} and \eqref{eq:IIa} with $\bm{w} = \bm{u}_n$ and $\rho = \psi_n$ and by also introducing the primitive $F$ of $f$, i.e. $F(r)=\int_0^r f(s) \;{\rm d}s$, we get
\begin{align*}
    &\frac{\d}{\d t} \Bigl[ \frac{1}{2} \lVert \Delta\phi_n\rVert_{L^2(Q)}^2 -\lVert \nabla\phi_n\rVert_{L^2(Q)}^2+ (F(\phi_n),1)+ \frac{1}{2M} \lVert \bm{u}_n\rVert_{\mathbb{H}}^2 \Bigr]\nonumber\\[2mm] &\,+\frac{1}{M}\lVert \sqrt{\eta(\phi_n)}D\bm{u}_n\rVert_{\mathbb{H}}^2 + \lVert \sqrt{m(\phi_n)}\nabla\psi_n\rVert_{L^2(Q)}^2=0.
\end{align*}
By using Assumption (A1), by noting that there exist positive constants $c_1$ and $c_2$ such that $F(s) \geqslant c_1 s^4 - c_2$, and by observing that
\begin{equation*}
-\lVert \nabla \phi_n \rVert_{L^2(Q)}^2 \geqslant -\frac{1}{4} \lVert \Delta \phi_n \rVert_{L^2(Q)}^2 - \frac{c_1}{2} \lVert \phi_n \rVert_{L^4(Q)}^4 - c,
\end{equation*}
we integrate over $(0,t)$, $t \leq T_n$, to get
\begin{align}
        &\frac{1}{4} \lVert \Delta \phi_n(t) \rVert_{L^2(Q)}^2 + \frac{c_1}{2} \lVert \phi_n(t) \rVert_{L^4(Q)}^4 + \frac{1}{2M} \lVert \bm{u}_n(t) \rVert_{\mathbb{H}}^2+ \frac{\eta_0}{M} \int_0^t \lVert D \bm{u}_n \rVert_{\mathbb{H}}\nonumber\\
        &\quad + m_0 \int_0^t \lVert \nabla \psi_n \rVert_{L^2(Q)}^2   \nonumber\\
        &\leqslant c + \frac{1}{2} \lVert \Delta \phi_{n}(0) \rVert_{L^2(Q)}^2 + (F(\phi_{n}(0)),1) + \frac{1}{2M} \lVert \bm{u}_{n}(0) \rVert_{\mathbb{H}}^2, \label{W18}
\end{align}
from which we first deduce that \eqref{W18} works also for $T_n = T$ and then we obtain the bound
\begin{equation}
\label{eq:Bphi}
    \lVert \bm{u}_n \rVert_{L^{\infty}(0,T;\mathbb{H})\cap L^2(0,T;\mathbb{V})} + \lVert \phi_n \rVert_{L^{\infty}(0,T;\Phi_2)}  \leqslant c.
\end{equation}
Here, we used Lemma \ref{poincphi} and \eqref{step0}. Moreover, from the estimate \eqref{W18} we get that $\nabla \psi_n$ is uniformly bounded in $L^2(0,T;L^2(Q))$. Hence, by using once more Lemma \ref{poincphi} together with \eqref{step00}, \eqref{W18}, and \eqref{eq:Bphi}, we find the bound
\begin{equation}
\label{eq:Bpsi}
\lVert \psi_n \rVert_{L^2(0,T;\Phi_1)} \leqslant c.
\end{equation}

$\bullet$ {\it Step 3:}
First, we note that equation \eqref{eq:Ia} can be written as follows
\begin{equation}
\label{314}
	\bm{u}_n' + P_{\mathbb{V}^n} \bigl( \eta(\phi_n) \bm{D}\bm{u}_n + B_0(\bm{u}_n,\bm{u}_n) + M \phi_n \nabla \psi_n		\bigr)=0 \quad \text{in} \; \mathbb{V}^n.
\end{equation}
Observe that estimates \eqref{eq:Bphi} and \eqref{eq:Bpsi} readily imply
\begin{equation}
\label{estimateR}
\lVert \phi_n \nabla \psi_n \rVert_{L^2(0,T;\mathbb{H})} \leq c .
\end{equation}
Moreover, by using \eqref{eq:Bphi} in \eqref{estimateB}, we get the bound on the inertia term, namely,
\begin{equation}
\label{estimateInertia}
\lVert B_0(\bm{u}_n,\bm{u}_n)\rVert_{L^{4/3}(0,T;\mathbb{V}')} \leqslant c,
\end{equation}
and by Assumption (A1), Korn inequality, and estimate \eqref{eq:Bphi}, we deduce
\begin{equation}
\label{estimateA}
\lVert \eta(\phi_n) \bm{D}\bm{u}_n \rVert_{L^2(0,T;\mathbb{H})} \leq c.
\end{equation}
Eventually, recalling that $\lVert P_{\mathbb{V}^n} \rVert_{\mathcal{L}(\mathbb{V}^n,\mathbb{V}^n)} \leq 1$, by comparison in equation \eqref{314}, we obtain the bound
\begin{equation}
\label{uprimeW}
\lVert \bm{u}_n' \rVert_{L^{4/3}(0,T;\mathbb{V}')}  \leqslant c.
\end{equation}
Similarly, we write equation \eqref{eq:IIa} as follows
\begin{equation}
\label{316}
	\phi_n' + P_{\Phi_2^n} \bigl( \bm{u}_n \cdot \nabla \phi_n - \nabla \cdot  ( m(\phi_n) \nabla \psi_n) \bigr) = 0 \quad \text{in} \; \Phi_2^n.
\end{equation}
By using the estimates \eqref{eq:Bphi} and \eqref{eq:Bpsi} it is straightforward to deduce
\begin{equation}
\label{estimateD}
\lVert \bm{u}_n \cdot \nabla \phi_n \rVert_{L^{\infty}(0,T;L^{6/5}(Q))} \leq c,
\end{equation}
as well as
\begin{equation}
\label{estimateC}
\lvert \nabla \cdot (m(\phi_n) \nabla \psi_n ) \rVert_{L^2(0,T;\Phi_1')} \leq c.
\end{equation}
Hence, by comparison in equation \eqref{316}, we get the bound
\begin{equation}
\label{phiprimeW}
\lVert \phi_n'\rVert_{L^2(0,T;\Phi_1')} \leqslant  c,
\end{equation}
where we used $\lVert P_{\Phi_2^n} \rVert_{\mathcal{L}(\Phi_2^n,\Phi_2^n)} \leq 1$.

$\bullet$ {\it Step 4:}
Take the scalar product in $L^2(Q)$ between equation \eqref{eq:IIIa} and $(\Delta^2 + 2\Delta)\phi_n$. We have
\begin{align*}
    \lVert \Delta^2 \phi_n \rVert_{L^2(Q)}^2 &+ 4 \lVert \Delta \phi_n \rVert_{L^2(Q)}^2 - 4\lVert \nabla \Delta \phi_n \rVert_{L^2(Q)}^2 \\&= - (f(\phi_n),(\Delta^2+2\Delta)\phi_n) + (\psi_n,(\Delta^2 + 2 \Delta) \phi_n).
\end{align*}
To estimate the right-hand side, we use the H\"older and the Young inequalities. We obtain
\begin{equation}
\label{89}
    \frac{1}{2} \lVert \Delta^2 \phi_n \rVert_{L^2(Q)}^2 + 3 \lVert \Delta \phi_n \rVert_{L^2(Q)}^2 - 4 \lVert \nabla \Delta \phi_n \rVert_{L^2(Q)}^2 \leqslant c \biggl(  \lVert f(\phi_n) \rVert_{L^2(Q)}^2 + \lVert \psi_n \rVert_{L^2(Q)}^2  \biggr).
\end{equation}
Concerning the third term on the left-hand side, we first integrate by parts and then we use the H\"older and the Young inequalities to get
\begin{equation*}
    -4 \lVert \nabla \Delta \phi_n \rVert_{L^2(Q)}^2 \geqslant  -\frac{1}{3} \lVert \Delta^2 \phi_n \rVert_{L^2(Q)}^2 - 12 \lVert \phi_n \rVert_{L^2(Q)}^2.
\end{equation*}
By using this result in \eqref{89}, we obtain
\begin{equation}
\label{deltaquadrofil2}
     \lVert \Delta^2 \phi_n \rVert_{L^2(Q)}^2 \leqslant c \biggl(  \lVert \Delta \phi_n \rVert_{L^2(Q)}^2  +  \lVert f(\phi_n) \rVert_{L^2(Q)}^2 + \lVert \psi_n \rVert_{L^2(Q)}^2  \biggr).
\end{equation}
Let us now consider the second term on the right-hand side in \eqref{deltaquadrofil2}. Thanks to \eqref{eq:Bphi}, we have
\begin{equation}
\label{boundfW}
    \lVert f(\phi_n) \rVert_{L^{\infty}(0,T;\Phi_2)} \leqslant c.
\end{equation}
Indeed, as $\phi$ is uniformly bounded in $\Phi_2$ and $f \in C^2(\mathbb{R})$ then $f(\phi_n)$ is uniformly bounded in $\Phi_2$. Thus, by using the estimates \eqref{eq:Bphi}, \eqref{eq:Bpsi}, and \eqref{boundfW} in \eqref{deltaquadrofil2}, we get
\begin{equation*}
    \lVert \Delta^2 \phi_n \rVert_{L^2(0,T;L^2(Q))} \leqslant c,
\end{equation*}
from which, by also making use once again of the estimate \eqref{eq:Bphi} and Lemma \ref{poincphi}, we deduce
\begin{equation}
\label{99}
    \lVert \phi_n \rVert_{L^2(0,T;\Phi_4)} \leqslant c.
\end{equation}
Eventually, by comparison in equation \eqref{eq:IIIa}, the estimates \eqref{eq:Bpsi}, \eqref{boundfW}, and \eqref{99} imply the bound
\begin{equation}
\label{106}
    \lVert \phi_n \rVert_{L^2(0,T;\Phi_5)} \leqslant c.
\end{equation}

\subsection{Passing to the limit}
\label{limitweak}
In this subsection we prove that the solution $(\bm{u}_n,\phi_n)$ of the finite-dimensional problem \eqref{eq:Ia}-\eqref{eq:IIIIa} converges to a weak solution of the problem \eqref{eq:ns}-\eqref{eq:periodic}.

Starting from estimates \eqref{eq:Bphi}, \eqref{eq:Bpsi}, \eqref{uprimeW}, \eqref{phiprimeW}, \eqref{106}, we find that, up to not relabeled subsequences, the following convergences hold
\begin{align}
&\bm{u}_n \overset{\ast}{\rightharpoonup} \bm{u}& &\text{in} \quad L^{\infty}(0,T;\mathbb{H}), \label{u1convW}\\
&\bm{u}_n \rightharpoonup \bm{u}& &\text{in} \quad L^2(0,T;\mathbb{V}),\label{u2convW}\\
&\bm{u}_n' \rightharpoonup \bm{u}'& &\text{in} \quad L^{4/3}(0,T;\mathbb{V}'),\label{uprimeconvW}\\
&\phi_n \overset{\ast}{\rightharpoonup} \phi& &\text{in} \quad L^{\infty}(0,T;\Phi_2),\label{phi1convW}\\
&\phi_n \rightharpoonup \phi& &\text{in} \quad L^2(0,T;\Phi_5),\label{phi2convW} \\
&\phi_n' \rightharpoonup \phi'& &\text{in} \quad L^2(0,T;\Phi_1'),\label{phiprimeconvW}\\
&\psi_n \rightharpoonup \psi& &\text{in} \quad L^2(0,T;\Phi_1).\label{psiconvW}
\end{align}
Furthermore, recalling, in particular, \eqref{uprimeW} and \eqref{phiprimeW}, the Aubin-Lions-Simon Lemma (see \cite[Thm.~3]{simon}) implies the following strong convergences
\begin{align}
    \bm{u}_n \rightarrow \bm{u}&  &\text{in}& \quad L^2(0,T;\mathbb{H}),\label{eq:ustrH}\\
    \bm{u}_n \rightarrow \bm{u}&  &\text{in}& \quad C^0([0,T];\mathbb{V}'),\label{eq:ustrH2}\\
    \phi_n   \rightarrow \phi&  &\text{in}& \quad C^0([0,T];\Phi_1).\label{eq:phistrH}
\end{align}
The weak convergences \eqref{u1convW}-\eqref {psiconvW} and the above strong convergences entail the following
\begin{align}
&(\bm{u}_n \cdot \nabla)\bm{u}_n \rightharpoonup (\bm{u}\cdot\nabla)\bm{u}& &\text{in}\quad L^{4/3}(0,T;\mathbb{L}^{6/5}(Q)),\label{convi}\\
&\eta(\phi_n)\bm{D}\bm{u}_n \rightharpoonup \eta(\phi) \bm{D}\bm{u}& &\text{in} \quad L^2(0,T;\mathbb{H}),\label{convii}\\
&\phi_n\nabla \psi_n \rightharpoonup \phi\nabla\psi& &\text{in} \quad L^2(0,T;\mathbb{H}),\label{conviii}\\
&\bm{u}_n\cdot \nabla \phi_n \overset{\ast}{\rightharpoonup} \bm{u}\cdot \nabla\phi& &\text{in} \quad L^{\infty}(0,T;L^{6/5}(Q)),\label{conviiv}\\
&m(\phi_n) \nabla \psi_n \rightharpoonup m(\phi)\nabla\psi& &\text{in} \quad L^2(0,T;L^2(Q)),\label{conviv}
\end{align}
Although convergences \eqref{convi}-\eqref{conviv} are obtained by standard arguments (see for instance \cite[Theorem VI.2.1, p.434]{Boyer}), we give here, for the sake of completeness, the idea of the proof.

$\bullet$ {\it Step 1: Strong convergence of nonlinear functions.}
The identification of the limit of the nonlinear functions $f$, $\eta$, and $m$, readily follows from the strong convergence \eqref{eq:phistrH}. Note indeed that $f$, $\eta$, and $m$ are Lipschitz continuous function and that the convergence \eqref{eq:phistrH} implies, by Sobolev embedding theorem, the strong convergence $\phi_n \rightarrow \phi$ in $L^{\infty}(0,T;L^6(Q))$. Then, recalling $\phi_n$ and $\phi$ are globally bounded (see \eqref{W18}), we have that
\begin{align*}
&\lVert f(\phi_n) - f(\phi) \rVert _{L^{\infty}(0,T;L^6(Q))} \leq c \lVert \phi_n - \phi \rVert_{L^{\infty}(0,T;L^6(Q))} \rightarrow 0,\\
&\lVert \eta(\phi_n) - \eta(\phi) \rVert _{L^{\infty}(0,T;L^6(Q))} \leq c \lVert \phi_n - \phi \rVert_{L^{\infty}(0,T;L^6(Q))} \rightarrow 0,\\
&\lVert m(\phi_n) - m(\phi) \rVert _{L^{\infty}(0,T;L^6(Q))} \leq c\lvert \phi_n - \phi \rVert_{L^{\infty}(0,T;L^6(Q))} \rightarrow 0,
\end{align*}
namely, $f(\phi_n) \rightarrow f(\phi)$ in $L^{\infty}(0,T;L^6(Q))$, $\eta(\phi_n) \rightarrow \eta(\phi)$ in $L^{\infty}(0,T;L^6(Q))$, and $m(\phi_n) \rightarrow m(\phi)$ in $L^{\infty}(0,T;L^6(Q))$, respectively.

$\bullet$ {\it Step 2: Weak convergence of the viscous term and \eqref{conviv}}
The strong convergences of $\eta(\phi_n)$ and $m(\phi_n)$ together with the weak convergences \eqref{u2convW} and \eqref{psiconvW} entail convergences \eqref{convii} and \eqref{conviv}. Indeed we deal, in both cases, with the product of a strongly converging sequence and a weakly converging one. Since this is sufficient to identify the limit functions (at least in $\mathcal{D}'((0,T)\times\Omega)$), we get convergences \eqref{convii} and \eqref{conviv} by taking advantage of the already proved estimates \eqref{estimateA} and \eqref{estimateC}, respectively. In particular, thanks to the estimates \eqref{estimateA} and \eqref{estimateC} we extract not relabeled converging subsequences $\eta(\phi_n)\bm{D}\bm{u}_n \rightharpoonup \bm{\zeta}$ in $L^2(0,T;\mathbb{H})$ and $m(\phi_n)\nabla\psi_n \rightharpoonup \varphi$ in $L^2(0,T;L^2(Q))$, respectively. Then, it is straightforward to show that, for any test function $\bm{v} \in L^2(0,T;\mathbb{W}^{1,\infty}(Q))$, one has
\begin{align*}
&\int_0^T ( \eta(\phi_n)\bm{D}\bm{u}_n - \eta(\phi)\bm{D}\bm{u},\nabla \bm{v} )\\
&= \int_0^T \int_Q \bigl( \eta(\phi_n) - \eta(\phi)\bigr) \bm{D}\bm{u}_{n} \cdot \nabla \bm{v} + \int_0^T \int_Q \eta(\phi) \bigl( \bm{D}\bm{u}_n - \bm{D}\bm{u}\bigr) \cdot \nabla \bm{v}\rightarrow 0,
\end{align*}
and that, for any test function $\rho \in L^2(0,T;W^{1,3}(Q))$, one gets
\begin{align*}
&\int_0^T ( m(\phi_n)\nabla\psi_n - m(\phi)\nabla\psi,\nabla \rho )\\
&= \int_0^T \int_Q \bigl( m(\phi_n) - m(\phi)\bigr) \nabla\psi_{n} \cdot \nabla \rho + \int_0^T \int_Q m(\phi) \bigl( \nabla\psi_n - \nabla\psi\bigr) \cdot \nabla \rho\rightarrow 0,
\end{align*}
namely, $\eta(\phi_n)\bm{D}\bm{u}_n \rightharpoonup \eta(\phi)\bm{D}\bm{u}$ in $L^2(0,T;(\mathbb{W}^{1,\infty}(Q))')$ and $m(\phi_n)\nabla\psi_n \rightharpoonup m(\phi)\nabla\psi$ in $L^2(0,T;(W^{1,3}(Q))')$, respectively. Then, by uniqueness of the limit, we obtain the identifications $\bm{\zeta} = \eta(\phi)\bm{D}\bm{u}$ and $\varphi = m(\phi)\nabla\psi$ which imply convergences \eqref{convii} and \eqref{conviv}, respectively.

$\bullet$ {\it Step 3: Weak convergence of the inertia term.}
Estimate \eqref{estimateInertia} entails the weak convergence (up to a not relabeled subsequence) $\bm{u}_n \cdot \nabla \bm{u}_n \rightharpoonup \bm{\xi}$ in $L^{4/3}(0,T;\mathbb{V}')$. First, via the strong convergence \eqref{eq:ustrH} of $\bm{u}_n$ and the weak convergence \eqref{u2convW} of $\nabla\bm{u}_n$, we deduce $\bm{\xi}=\bm{u}\cdot\nabla\bm{u}$ by arguing as above. Then, we refine this convergence by proving sharp estimates. We have
$$
\lVert \bm{u}_n \cdot \nabla \bm{u}_n \rVert_{L^2(0,T;\mathbb{L}^1(Q))} \leq \lVert \bm{u}_n \rVert_{L^{\infty}(0,T;\mathbb{H})} \lVert \bm{u}_n \rVert_{L^2(0,T;\mathbb{V})} \leq c,
$$
and
$$
\lVert \bm{u}_n \cdot \nabla \bm{u}_n \rVert_{L^1(0,T;\mathbb{L}^{3/2}(Q))} \leq \lVert \bm{u}_n \rVert_{L^2(0,T;\mathbb{V})}^2 \leq c.
$$
Hence, by standard interpolation theory, as
$$
(L^2(0,T;\mathbb{L}^1(Q)),L^1(0,T;\mathbb{L}^{3/2}(Q)))_{1/2} \subset L^{4/3}(0,T;\mathbb{L}^{6/5}(\Omega)),
$$
we eventually obtain \eqref{convi}.

$\bullet$ {\it Step 4: Convergence of Korteweg force and advective term.}
Convergences \eqref{conviii} and \eqref{conviiv} follows from convergences \eqref{psiconvW}, \eqref{eq:ustrH}, and \eqref{eq:phistrH} and estimates \eqref{estimateR} and \eqref{estimateD} by arguing as in {\it Step 2}.

We now prove that the limit functions $\bm{u}$, $\phi$, and $\psi$ satisfy \eqref{deboleteorema1}-\eqref{deboleteorema3}. To this aim we pass to the limit in the weak formulation \eqref{eq:Ia}-\eqref{eq:IIIa}. For any $k$ fixed, let $\bm{v}_k \in \mathbb{V}^k$, $\rho_k \in \Phi_2^k$, and $\theta \in \mathcal{D}(0,T)$. We have
\begin{equation*}
\begin{split}
&\int_0^T (\bm{u}_n',\bm{v}_k) \theta(t) + \int_0^T (\eta(\phi_n)\bm{D}\bm{u}_n,\nabla\bm{v}_k) \theta(t)+ \int_0^T (B_0(\bm{u}_n,\bm{u}_n),\bm{v}_k) \theta(t) \\&\qquad+ M \int_0^T (\phi_n \nabla \psi_n,\bm{v}_k) \theta(t) =0,\\
& \int_0^T (\phi_n',\rho_k)\theta(t) + \int_0^T (\bm{u}_n \cdot \nabla \phi_n,\rho_k)\theta(t) + \int_0^T (m(\phi_n)\nabla\psi_n,\nabla\rho_k)\theta(t) =0,\\[2mm]
& (\Delta^2\phi_n,\rho_k) + 2(\Delta\phi_n,\rho_k)+(f(\phi_n),\rho_k) - (\psi_n,\rho_k)=0 \qquad \text{a.e.}\;(0,T).
\end{split}
\end{equation*}
Since $k$ is fixed we pass to the limit $n \rightarrow \infty$ and, by using weak convergences proved above, we obtain
\begin{equation*}
    \begin{split}
        &\int_0^T (\bm{u}',\bm{v}_k)\theta(t) + \int_0^T (\eta(\phi)\bm{D}\bm{u},\nabla\bm{v}_k)\theta(t) + \int_0^T (B_0(\bm{u},\bm{u}),\bm{v}_k)\theta(t) \\&\qquad+ M \int_0^T (\phi \nabla \psi,\bm{v}_k)\theta(t)=0, \\
        &\int_0^T (\phi',\rho_k)\theta(t) + \int_0^T (\bm{u}\cdot \nabla \phi,\rho_k)\theta(t) + \int_0^T (m(\phi)\nabla\psi,\nabla \rho_k)\theta(t)=0, \\[2mm]
        &(\Delta^2\phi,\rho_k) + 2(\Delta\phi,\rho_k)+(f(\phi),\rho_k) - (\psi,\rho_k)=0 \qquad\text{a.e. in}\;(0,T).
    \end{split}
\end{equation*}
Now let $\bm{v} \in \mathbb{V}$, $\rho \in \Phi_2(Q)$ and set $\bm{v}_k=P_{\mathbb{V}^k}\bm{v}$, $\rho_k=P_{\Phi_2^k}\rho$. As $\bm{v}_k \rightarrow \bm{v}$ in $\mathbb{V}$ and $\rho_k \rightarrow \rho$ in $\Phi_2$, we pass to the limit $k\rightarrow \infty$ and we obtain
\begin{align}
        &\int_0^T (\bm{u}',\bm{v})\theta(t) + \int_0^T (\eta(\phi)\bm{D}\bm{u},\nabla\bm{v})\theta(t) + \int_0^T (B_0(\bm{u},\bm{u}),\bm{v})\theta(t)\nonumber \\&\qquad+ M \int_0^T (\phi\nabla\psi,\bm{v})\theta(t)=0, \label{equationtheta}\\
        &\int_0^T (\phi',\rho)\theta(t) + \int_0^T (\bm{u}\cdot\nabla\phi,\rho)\theta(t) + \int_0^T (m(\phi)\nabla\psi,\nabla\rho)\theta(t)=0,\label{equationtheta0} \\[2mm]
        &(\Delta^2\phi,\rho) + 2(\Delta\phi,\rho)+(f(\phi),\rho) - (\psi,\rho)=0 \qquad\text{a.e. in}\;(0,T),\label{equationtheta1}
\end{align}
where, in particular, equation \eqref{equationtheta1} implies \eqref{deboleteorema3}. Then, as the set of functions of the form $\bm{v}\theta$ is dense in $ \mathcal{D}(0,T;\mathbb{V})$ and the set of functions of the form $\rho\theta$ is dense in $\mathcal{D}(0,T;\Phi_2)$, from equations \eqref{equationtheta} and \eqref{equationtheta0}, we deduce that $\bm{u}$, $\phi$, and $\psi$ satisfy \eqref{deboleteorema1} and \eqref{deboleteorema2} almost everywhere in $(0,T)$.

{\color{black}
Eventually, we prove the convergence of the initial conditions \eqref{deboleteorema4}. The strong convergence \eqref{eq:ustrH2} implies that $\bm{u}_n(0) \rightarrow \bm{u}(0)$ in $\mathbb{V}'$. On the other hand, by definition of $\mathbb{V}^n$, we have the convergence $\bm{u}_n(0)=P_{\mathbb{V}^n}\bm{u}_0 \rightarrow \bm{u}_0$ in $\mathbb{H}$. As $\mathbb{H} \subset \mathbb{V}'$, by uniqueness of the limit in $\mathbb{V}'$, we deduce $\bm{u}(0) = \bm{u}_0$. Similarly, from the strong convergence \eqref{eq:phistrH}, we argue as above to obtain $\phi(0) = \phi_0$.
}

Finally, to establish the energy inequality \eqref{enineq} we can take $\bm{v} = \bm{u}_{n} $  in \eqref{eq:Ia} and $\rho = \psi_{n}$ in \eqref{eq:IIa} to get
\begin{align*}
    \frac{1}{2M}\|\bm{u}_{n}(t)\|_{\mathbb{H}}^2 &+ \mathcal{F}_{\text{sh}}[\phi_{n}(t)]  + \frac{1}{M}\int_0^t \lVert \sqrt{\eta(\phi_{n})} D \bm{u}_{n} \rVert_{\mathbb{H}}^{2} +
 \int_0^t \lVert \sqrt{m(\phi_{n})} \nabla \psi_{n} \rVert_{L^2(Q)}^2   \\ &=  \frac{1}{2}\|\bm{u}_{n}(0)\|_{\mathbb{H}}^2 + \mathcal{F}_{\text{sh}}[P_{\Phi_2^n} \phi_0]\\
    &\le \frac{1}{2}\|\bm{u}_0\|_{\mathbb{H}}^2 + \mathcal{F}_{\text{sh}}[P_{\Phi_2^n} \phi_0].
\end{align*}
Then, we follow, for instance, the argument detailed in \cite[Chapt.V, Proof of Prop.V.1.7]{Boyer}, using also the fact that
$\mathcal{F}_{\text{sh}}[P_{\Phi_2^n} \phi_0] \to \mathcal{F}_{\text{sh}}[\phi_0]$
as $n$ goes to $\infty$.


\subsection{Dimension two: uniqueness (sketch proof)}
\label{weakunique2D}

We now want to prove that weak solutions are unique in the case $d=2$. We begin by taking an appropriate test function in the first equation of the weak formulation satisfied by the difference of two arbitrary weak solutions. Simplifying the resulting equation will lead to an inequality of the form (see \eqref{d1} below)

\begin{equation*}
\frac{1}{2}\|\bm{u}(t)\|_{\mathbb{L}^2(Q)}^{2} + c\int_0^t \|\nabla \bm{u}(s)\|_{\mathbb{L}^2(Q)}^{2} ~\mathrm{d}s \le \frac{1}{2}\|\bm{u}(0)\|_{\mathbb{L}^2(Q)}^{2} + \sum_{j=1}^{8}I_{j}.
\end{equation*}
Next, we bound each of the terms $I_{j}$, $j=1,\dots,8$. We approach the second equation of the weak formulation in a similar way, leading to (see \eqref{d3} below)

\begin{equation*}
    \frac{1}{2}\|\nabla \phi(t)\|_{\mathbb{L}^2(Q)}^{2} + \int_0^t \|\Delta^{2}\phi(s)\|_{\mathbb{L}^2(Q)}^{2}~\mathrm{d}s = \frac{1}{2}\|\nabla \phi(0)\|_{\mathbb{L}^2(Q)}^{2} + \sum_{j=9}^{15}I_{j}.
\end{equation*} Estimating $I_{1}, \dots,I_{15}$ and adding up the resulting inequalities eventually allows us to apply Osgood's lemma (\cite{JFAGiorg}, Appendix B) which will complete the proof. Note that in the following of this section we adopt the notation $\| \cdot \| \equiv \| \cdot \|_{\mathbb{L}^{2}(Q)}$ for the sake  of simplicity.

\subsection{The momentum equation}
 We now begin working in accordance with the above plan. Let $(\bm{u}_{1}, \phi_{1}), ~(\bm{u}_{2}, \phi_{2})$ be two weak solutions, and define $(\bm{u}, \phi) := (\bm{u}_{1} - \bm{u}_{2}, \phi_{1} - \phi_{2})$ to be the difference of the two solutions. 

 The first equation in our weak formulation  is equivalent to
\begin{gather}
(\bm{u}', \bm{v}) + b_{0}(\bm{u}, \bm{u}, \bm{v}) + ( \eta(\phi)D(\bm{u}), \nabla \bm{v}) = - M(\phi \nabla \psi , \bm{v}) \hspace{10pt} \forall \hspace{3pt} \bm{v} \in L^{2}(0,T; \mathbb{V}), \vspace{5pt} \label{x1}
\end{gather} almost everywhere in $(0,T)$. This is satisfied by $(u_{i}, \phi_{i})$ for $i=1,2$.
{\color{black} Subtracting  \eqref{x1} with $i=2$ from \eqref{x1} with $i=1$ gives us
\begin{equation*} \begin{aligned}
     (\bm{u}',\bm{v}) &+ b_{0}(\bm{u}_{1},\bm{u}_{1},\bm{v}) - b_{0}(\bm{u}_{2},\bm{u}_{2},\bm{v}) + (\eta(\phi_{1})D(\bm{u}_{1}), \nabla \bm{v}) - (\eta(\phi_{2})D(\bm{u}_{2}), \nabla \bm{v}) \\[1ex] & + M(\phi_{1}\nabla \psi_{1}, \bm{v}) - M(\phi_{2}\nabla \psi_{2}, \bm{v}) = 0.
\end{aligned}
\end{equation*}
Taking then $\bm{v} = \bm{u}$ and integrating in time,
\begin{equation}
\begin{aligned}
    \frac{1}{2}\|\bm{u}(t)\|^{2} &+ \int_0^t  b_{0}(\bm{u}_{1}, \bm{u}_{1}, \bm{u}) - \int_0^t  b_{0}(\bm{u}_{2}, \bm{u}_{2}, \bm{u}) \notag + \int_0^t  (\eta(\phi_{1})D(\bm{u}_{1}), \nabla \bm{u}) \\
    \quad & - \int_0^t  (\eta(\phi_{2})D(\bm{u}_{2}), \nabla \bm{u}) \notag + M\int_0^t  (\phi_{1} \nabla \psi_{1}, \bm{u})  - M\int_0^t  (\phi_{2} \nabla \psi_{2}, \bm{u}) = \frac{1}{2}\|\bm{u}(0)\|^{2}.
\end{aligned}
\end{equation}
}Next, we rewrite each of the integrals appearing above into a form which will allow us to carry out effective estimates. Using the trilinearity of the form $b_{0}$, we have that  $b_{0}(\bm{u}_{1}, \bm{u}_{1}, \bm{u}) - b_{0}(\bm{u}_{2}, \bm{u}_{2}, \bm{u}) =  b_{0}(\bm{u}, \bm{u}_{1}, \bm{u})$. Additionally, we note that by adding and subtracting the quantity $\eta(\phi_{1})D(\bm{u}_{2}(t))$, the following equality holds  \begin{align*}
 (\eta(\phi_{1})D(\bm{u}_{1}(t)), \nabla \bm{u}(t)) &- (\eta(\phi_{2})D(\bm{u}_{2}(t)), \nabla \bm{u}(t)) \\[3mm] =  &~(\eta(\phi_{1})D(\bm{u}(t)), D(\bm{u}(t)))+ ((\eta(\phi_{1}) - \eta(\phi_{2}))D(\bm{u}_{2}(t)), \nabla \bm{u}(t)),
\end{align*}
which will be useful for the subsequent analysis. Observe now that the Korteweg force term $\phi \nabla \psi$ can be re-expressed as follows
\begin{align*}
    \phi \nabla \psi &= \nabla (\phi \psi) - \psi \nabla \phi = \nabla (\phi \psi) - (f(\phi) + 2\Delta \phi + \Delta^{2} \phi)\nabla \phi \\
    &= \nabla (\phi \psi) - \nabla  F(\phi) - 2 \Delta \phi \nabla \phi - \Delta^{2}\phi \nabla \phi.
\end{align*} In component notation, we have
\begin{equation*}
    (\Delta^{2}\phi \nabla \phi)_{i} = \text{div}(\nabla \Delta \phi)\partial_{i}\phi = \partial_{j}(\partial_{j}\Delta \phi)\partial_{i}\phi = \partial_{j}(\partial_{j} \Delta \phi \partial_{i}\phi) - \partial_{j}\Delta \phi \hspace{2pt} \partial_{i} \partial_{j} \phi,
\end{equation*} so that
    $\Delta^{2}\phi \nabla \phi = \text{div}(\nabla \phi \otimes \nabla \Delta \phi) - (\nabla \Delta \phi \cdot \nabla) \nabla \phi$,
 where $\mathbb{R}^{d \times d} \ni (v \otimes w)_{ij} := v_{i}w_{j}$, for $v, w \in \mathbb{R}^{d}$. Thus, we get
\begin{align*}
    (\phi_{1}\nabla \psi_{1}, \bm{u}) &= \int_{Q} \left(\nabla (\phi_{1} \psi_{1}) - \nabla F(\phi) \right) \bm{u}  \\ \quad  & - 2 \int_{Q} \Delta \phi_{1} \nabla \phi_{1} \cdot \bm{u}  - \int_{Q} \Delta^{2} \phi_{1} \nabla \phi_{1} \cdot \bm{u}.
\end{align*}
It follows from integration by parts and the divergence free condition imposed on $\bm{u}$ that the first integral on the right hand-side of the above equality is zero. Using our previous observations and integration by parts, we also have
\begin{align*}
    &\int_{Q} \Delta^{2} \phi_{1} \nabla \phi_{1} \cdot \bm{u}   = - \int_{Q} (\nabla \phi_{1} \otimes \nabla \Delta \phi_{1}) : \nabla \bm{u} - b_{0}(\nabla \Delta \phi_{1}, \nabla \phi_{1}, \bm{u}), \\
    &\int_{Q} \Delta \phi_{1} \nabla \phi_{1} \cdot \bm{u} = - \int_{Q} (\nabla \Delta \phi_{1} \cdot \bm{u}) \phi_{1} ,
\end{align*} and similarly for $\phi_{2}$. Therefore, we find that
\begin{align*}
    &(\phi_{1} \nabla \psi_{1}, \bm{u}) - (\phi_{2} \nabla \psi_{2}, \bm{u}) = -\int_{Q} \nabla \phi_{2} \otimes \nabla \Delta \phi_{2} : \nabla \bm{u}  + \int_{Q} \nabla \phi_{1} \otimes \nabla \Delta \phi_{1} : \nabla \bm{u} \\ &~~+ b_{0}(\nabla \Delta \phi_{2}, \nabla \phi_{2}, \bm{u}) - b_{0}(\nabla \Delta \phi_{1}, \nabla \phi_{1}, \bm{u}) - 2 \int_{Q} (\nabla \Delta \phi_{2} \cdot \bm{u})\phi_{2} + 2 \int_{Q} (\nabla \Delta \phi_{1} \cdot \bm{u})\phi_{1} .
\end{align*}
Exploiting linearity we can rewrite this equality as
\begin{align*}
    (\phi_{1} \nabla \psi_{1}&, \bm{u}) - (\phi_{2} \nabla \psi_{2}, \bm{u}) = - \int_{Q} \nabla \phi_{2} \otimes \nabla \Delta \phi : \nabla \bm{u} - \int_{Q} \nabla \phi \otimes \nabla \Delta \phi_{1} : \nabla \bm{u}  \\ &- b_{0}(\nabla \Delta \phi, \nabla \phi_{2}, \bm{u}) - b_{0}(\nabla \Delta \phi_{1}, \nabla \phi, \bm{u}) - \int_{Q} (\nabla \Delta \phi \cdot \bm{u}) \phi_{2}  - \int_{Q} (\nabla \Delta \phi_{1} \cdot \bm{u}) \phi.
\end{align*}
Thus, we obtain the identity
\begin{equation}
    \begin{aligned}
         &\frac{1}{2}\|\bm{u}(t)\|^{2} + \int_0^t   \eta(\phi_{1})\|D\bm{u}\|^{2}  = \frac{1}{2}\|\bm{u}(0)\|^{2} + \int_0^t  b_{0}(\bm{u}, \bm{u}_{1}, \bm{u})  \notag \\ &+ \int_0^t  (\eta(\phi_{1}) - \eta(\phi_{2}))D(\bm{u}_{2}(t)), \nabla \bm{u}(t))  + M \int_0^t  (\nabla \phi_{2} \otimes \nabla \Delta \phi, \nabla \bm{u})   \\
    &+ M\int_0^t  (\nabla \phi \otimes \nabla \Delta \phi_{1}, \nabla \bm{u})   + M\int_0^t  b_{0}(\nabla \Delta \phi, \nabla \phi_{2}, \bm{u})  + M\int_0^t  b_{0}(\nabla \Delta \phi_{1}, \nabla \phi, \bm{u})   \\
    & + M\int_0^t  (\nabla \Delta \phi \cdot \bm{u}, \phi_{2})  + M\int_0^t  (\nabla \Delta \phi_{1} \cdot \bm{u}, \phi)
    =  \frac{1}{2}\|\bm{u}(0)\|^{2} + \sum_{j=1}^{8}I_{j}.
    \end{aligned}
\end{equation} This implies that
\begin{align} \label{d1}
    &\frac{1}{2}\|\bm{u}(t)\|^{2} + \frac{\eta_{0}}{\sqrt{2}} \int_0^t  \|\nabla \bm{u}\|^{2}  \le \frac{1}{2}\|\bm{u}(0)\|^{2} + \sum_{j=1}^{8}I_{j},
\end{align}
where we have used Korn's inequality $\|\nabla \bm{u}\| \le \sqrt{2} \|\Delta \bm{u}\|$ for $\bm{u} \in \mathbb{V}$ and the assumption on the viscosity term $0 < \eta_{0} < \eta(s)$ for each $s \in \mathbb{R}$. We now proceed to estimate $I_{1},\dots , I_{8}$. In what follows we will repeatedly denote by $c$ a positive constant possibly dependent on $\eta_{i}, m_{i}, M_0, M$ for $i \in \{1,2,3\}$, while $\epsilon, \epsilon '$ will be used to denote positive constants which arise from Young's inequality and can be chosen as small as required. We will also denote by $R$ a time-dependent function $R \in C([0,T])$ which may also depend on $\eta_{i}, m_{i}, M_0, M$ for $i \in \{1,2,3\}$. Let us also observe that since $\overline{\phi} = \overline{\phi_{1}(t)} - \overline{\phi_{2}(t)} = \overline{\phi_{0}} - \overline{\phi_{0}} =0$, we have that
$$
\|\phi\|_{H^{2}} \le c\|\Delta \phi\|, \quad\|\phi\|_{H^{3}} \le c\|\nabla \Delta \phi\|, \quad \|\phi\|_{H^{4}} \le c \|\Delta^{2}\phi\|, \quad \|\phi\|_{H^{1}} \le c\|\nabla \phi\|.
$$
We can now start estimating the terms $I_j$, $j=1,\dots,8$.

\noindent
$\bullet$	 $I_1$ $\bullet$
Using the H\"{o}lder, the Ladyzhenskaya, and the Young inequalities, we have
\begin{align*}
 {I}_{1} &\le  \int_0^t  \|\bm{u}\|_{L^{4}(Q)}^{2}  \| \nabla \bm{u}_{1}\| \le \frac{c}{4\epsilon} \int_0^t  \|\bm{u}\|^{2}\|\nabla \bm{u}_{1}\|^{2} + \epsilon \int_0^t  \|\bm{u}\|^{2}_{H^{1}(Q)}.
 \end{align*} Poincar\'{e}'s inequality implies that $\|\bm{u}\|^{2}_{H^{1}(Q)} \le C\|\nabla \bm{u}\|^{2}$. Thus,
 \begin{align*}
     I_{1} \le  c \epsilon \int_0^t  \|\nabla \bm{u}\|^{2}  + c\int_0^t  \|\bm{u}\|^{2}\|\nabla \bm{u}_{1}\|^{2} .
 \end{align*}
$\bullet$	 $I_2$ $\bullet$ By applying H\"{o}lder's inequality, we find
\begin{align*}
 I_{2} &\le \int_0^t  \|\eta(\phi_{1})- \eta(\phi_{2})\|_{L^{\infty}(Q)}\|D\bm{u}_{2}\| \|\nabla \bm{u}\|.
 \end{align*} On the other hand, we have \begin{equation*}
     \|\eta(\phi_{1})-\eta(\phi_{2})\|_{L^{\infty}(Q)} = \left\|\int_{0}^{1} \eta'(s\phi_{1} + (1-s)\phi_{2})\phi(t) \hspace{2pt} \mathrm{d}s \right\|_{L^{\infty}(Q)} \le C\|\phi\|_{L^\infty (Q)}.
 \end{equation*} Then, additionally, using the Br\'{e}zis-Gallouet and the Young inequalities, we get
 \begin{equation*}
     I_{2} \le \epsilon \int_0^t  \|\nabla \bm{u}\|^{2} + c\int_0^t  \|\nabla \phi\|^{2} \ln\left(\frac{c}{\|\nabla \phi\|^{2}} \right) \|D\bm{u}_{2}\|^{2}.
 \end{equation*} \\[2mm]
$\bullet$	 $I_3$ $\bullet$
Applying the H\"{o}lder and the Young inequalities yields
\begin{align}
 I_{3} &= M\int_0^t  (\nabla \phi_{2} \otimes \nabla \Delta \phi, \nabla \bm{u}) \le \epsilon \int_0^t  \|\nabla \bm{u}\|^{2} + \frac{1}{4\epsilon M^{2}} \int_0^t  \|\nabla \phi_{2}\|^{2}_{L^{\infty}(Q)} \|\nabla \Delta \phi\|^{2}.
 \label{I3a0}
 \end{align}  Firstly, we observe that,  integration by parts, we have
 \begin{equation} \label{I3a}
     \|\nabla \Delta \phi\|^{2} = (\nabla \Delta \phi, \nabla \Delta \phi) = - (\Delta \phi, \Delta^{2} \phi) \le \|\Delta \phi\|\|\Delta^{2} \phi\|.
 \end{equation} Therefore, using \eqref{I3a}, \vspace{-2mm} \begin{equation}
 \label{I3abis}
     \|\nabla \phi_{2}\|^{2}_{L^{\infty}(Q)}\|\nabla \Delta \phi\|^{2} \le \|\nabla \phi_{2}\|^{2}_{L^{\infty}} \|\Delta \phi\| \|\Delta^{2} \phi\| \le \epsilon' \|\Delta^{2} \phi\|^{2} + \frac{1}{4\epsilon'}\|\nabla \phi_{2}\|^{4}_{L^{\infty}} \|\Delta \phi\|^{2}.
 \end{equation} Also, using integration by parts once more, we get
 \begin{equation*}
     \|\Delta \phi\|^{2} = (\Delta \phi, \Delta \phi) = - (\nabla \phi, \nabla \Delta \phi) \le \|\nabla \phi\| \| \nabla \Delta \phi\| \le \|\nabla \phi\| \| \Delta \phi\|^{\frac{1}{2}}\|\Delta^{2} \phi\|^{\frac{1}{2}},
 \end{equation*} which implies, by interpolation, \begin{equation} \label{I3c}
     \|\Delta \phi\|^{2} \le \|\nabla \phi\|^{\frac{4}{3}}\|\Delta^{2} \phi\|^{\frac{2}{3}}.
 \end{equation} In light of this estimate, from \eqref{I3abis} we deduce the following
 \begin{align*}
     \|\nabla \phi_{2}\|^{2}_{L^{\infty}}\|\nabla \Delta \phi\|^{2} &\le \epsilon'\|\Delta^{2}\phi\|^{2} + \frac{1}{4\epsilon'}\|\nabla \phi_{2}\|^{4}_{L^{\infty}} \|\nabla \phi\|^{\frac{4}{3}}\|\Delta^{2}\phi\|^{\frac{2}{3}}
 \end{align*} Applying Young's inequality with $p=3/2, ~q=3$ yields \begin{equation}
 \label{I3b}
      \|\nabla \phi_{2}\|^{2}_{L^{\infty}}\|\nabla \Delta \phi\|^{2} \le 2\epsilon' \| \Delta^{2} \phi\|^{2} + C\|\nabla \phi_{2}\|^{6}_{L^{\infty}} \|\nabla \phi\|^{2}
 \end{equation}
 Combining \eqref{I3a0} with \eqref{I3b} and taking $\epsilon = \frac{\eta_{0}}{2\sqrt{2}}$, gives
 \begin{align*}
     I_{3} \le \frac{\eta_{0}}{2\sqrt{2}}\int_0^t  \|\nabla \bm{u}\|^{2} + \frac{\epsilon \sqrt{2}}{\eta_{0}} \int_0^t  \| \Delta^{2} \phi \|^{2}  + c \int_0^t  \|\nabla \phi_{2}\|^{6}_{L^{\infty}(Q)} \|\nabla \phi\|^{2}.
 \end{align*}

Note that, from \eqref{I3c}, we deduce the following two facts which will be useful for later estimates \begin{align}
     &\|\Delta \phi\|^{2} \le \|\nabla \phi\|^{2} + \|\Delta^{2}\phi\|^{2}, \label{Del}  \\[2mm]
     & \|\nabla \Delta \phi \|^{2} \le 2 \epsilon '\| \Delta^{2} \phi\|^{2} + C\|\nabla \phi\|^{2}. \label{Nab}
 \end{align}

$\bullet$	 $I_4$ $\bullet$
We proceed in a similar fashion as for $I_{3}$. Using H\"{o}lder, Ladyzhenskaya and Young's inequalities, we find
\begin{align*}
 I_{4} &= M\int_0^t  (\nabla \phi \otimes \nabla \Delta \phi_{1}, \nabla \bm{u})  \le c\int_0^t  \|\nabla \Delta \phi_{1}\|_{L^{4}(Q)}\|\nabla \phi\|_{L^{4}(Q)}\|\nabla \bm{u}\|  \\ \quad& \le \epsilon \int_0^t  \|\nabla \bm{u}\|^{2}   + c\int_0^t  \|\nabla \Delta \phi_{1}\|^{2}_{L^{4}(Q)}\|\Delta \phi\|^{2}.
 \end{align*} We now need to estimate $\| \Delta \phi\|^{2}$. Recalling the estimate $ \|\Delta \phi\|^{2} \le \|\nabla \phi\| \| \nabla \Delta \phi\|$, using Young's inequality and \eqref{Nab}, we get
 \begin{align*}
     \|\nabla \Delta \phi_{1}\|^{2}_{L^{4}(Q)} \|\Delta \phi\|^{2} \le 2 \epsilon' \|\Delta^{2}\phi\|^{2} + C\|\nabla \Delta \phi_{1}\|^{4}_{L^{4}(Q)} \|\nabla \phi\|^{2}.
 \end{align*}
Thus we arrive at
\begin{equation*}
     I_{4} \le \epsilon \int_0^t  \|\nabla \bm{u}\|^{2} + c \epsilon' \int_0^T  \|\Delta^{2} \phi\|^{2} + c \int_0^t  \|\nabla \Delta \phi_{1}\|^{4}_{L^{4}(Q)} \|\nabla \phi\|^{2}.
\end{equation*}
$\bullet$	 $I_5$ $\bullet$
Since $\phi_{2} \in L^{2}(0,T; H^{4}_{N}(Q))$, using H\"{o}lder's inequality, we obtain
\begin{align*}
 I_{5} = M\int_{0}^{t} b(\nabla \Delta \phi, \nabla \phi_{2}, \bm{u})  &\le c\int_{0}^{t} \|\nabla \Delta \phi\|_{L^{4}(Q)} \|\nabla^{2}\phi_{2}\| \|\bm{u}\|_{L^{4}(Q)} \\[1ex] &\le c \int_{0}^{t} \|\nabla \Delta \phi\|_{L^{4}(Q)} \|\bm{u}\|_{L^{4}(Q)} .
 \end{align*} Furthermore, note that there exists $c>0$ such that $\|\nabla \Delta \phi\|_{H^{1}(Q)} \le \|\phi\|_{H^{4}(Q)} \le c\|\Delta^{2}\phi\|$. This is a consequence of the fact that $\overline{\phi} = 0$ and the Poincar\'{e}-Wirtinger inequality. Thus, by Young's inequality,
 \begin{align*}
     I_{5} \le c \int_{0}^{t} \|\nabla \Delta \phi\|^{\frac{1}{2}}\|\Delta^{2} \phi\|^{\frac{1}{2}}\|\bm{u}\|^{\frac{1}{2}}\|\nabla \bm{u}\|^{\frac{1}{2}}  \le \epsilon \int_{0}^{t} \|\nabla \bm{u}\|^{2}  + \frac{1}{4\epsilon}\int_{0}^{t} \|\nabla \Delta \phi\|\|\Delta^{2}\phi\| ,
 \end{align*} where we have used the regular Poincar\'{e} inequality to deduce $\|\bm{u}\|\|\nabla \bm{u}\| \le \|\nabla \bm{u}\|^{2}$. Using Young's inequality and \eqref{Nab} yields
 \begin{align*}
     I_{5} \le \epsilon \int_{0}^{t} \|\nabla \bm{u}\|^{2}  + c \int_{0}^{t} \|\Delta^{2}\phi\|^{2} + c \int_{0}^{t} \|\nabla \phi\|^{2} .
 \end{align*}
$\bullet$ $I_6$ $\bullet$
Proceeding in a similar way as for $I_{5}$, H\"{o}lder's inequality gives
\begin{align*}
 I_{6}  \le c\int_{0}^{t} \|\nabla \Delta \phi_{1}\|_{L^{4}(Q)} \|\nabla\phi\|_{H^1(Q)}\|\bm{u}\|_{L^{4}(Q)}  \le c\int_{0}^{t} \|\nabla \Delta \phi_{1}\|_{L^{4}(Q)}\|\Delta \phi\| \|\bm{u}\|^{\frac{1}{2}}\|\bm{u}\|^{\frac{1}{2}}_{H^{1}(Q)},
 \end{align*} where we have used $\|\nabla \phi\|_{H^1(Q)} \le \|\phi\|_{H^{2}(Q)} \le C\|\Delta \phi\|$. Next, recalling \eqref{I3c} and using Poincar\'{e}'s inequality and Young's inequality with $p = 4, ~q=4/3$, we get
 \begin{align*}
     I_{6} &\le  \frac{\epsilon^{4}}{4} \int_{0}^{t} \|\nabla \bm{u}\|^{2}  + c\int_{0}^{t} \| \nabla \Delta \phi_{1}\|^{\frac{4}{3}}_{L^{4}(Q)} \|\nabla \phi\|^{\frac{8}{9}}\|\Delta^{2}\phi\|^{\frac{4}{9}}\|\bm{u}\|^{\frac{2}{3}}.
 \end{align*}
 Applying Young's inequality two more times (first with $p=9/2, ~q=9/7$ then with $p=7/4, ~q=7/3$), yields
 \begin{equation*}
     I_{6} \le \frac{\epsilon^{4}}{4} \int_{0}^{t} \|\nabla \bm{u}\|^{2} + c \int_{0}^{t} \|\Delta^{2}\phi\|^{2} + c \int_{0}^{t} \|\nabla \Delta \phi_{1}\|^{\frac{12}{7}}_{L^{4}} (\|\nabla \phi\|^{2} + \|\bm{u}\|^{2}) .
 \end{equation*}
$\bullet$  $I_7$ $\bullet$
Using the fact that $\phi_{2} \in L^{\infty}(0,T; H^{2}(Q)) \hookrightarrow L^{\infty}(0,T; L^{\infty}(Q))$ as well as the H\"{o}lder and Young inequalities, we get
\begin{align*}
 I_{7} = \int_{0}^{t} (\nabla \Delta \phi \cdot \bm{u}, \phi_{2}) \le \epsilon \int_{0}^{t} \|\Delta^{2}\phi\|^{2} + c \int_{0}^{t} \|\nabla \phi\|^{2}  + c \int_{0}^{t} \|\bm{u}\|^{2} .
\end{align*} \\[1mm]
$\bullet$	 $I_8$ $\bullet$
An application of H\"{o}lder and Young's inequalities entails that
\begin{align*}
 I_{8} \le c\int_{0}^{t} \|\nabla \Delta \phi_{1}\|_{L^{4}(Q)} \|\bm{u}\| \|\phi\|_{L^{4}(Q)}   \le c \int_{0}^{t} \|\nabla \Delta \phi_{1}\|^{2}_{L^{4}} \|\bm{u}\|^{2}   + \epsilon \int_{0}^{t} \|\phi\|^{2}_{L^{4}(Q)}.
 \end{align*} The Ladyzhenskaya and Poincar\'{e}-Wirtinger inequalities imply that \begin{equation} \label{I8a}
     \|\phi\|^{2}_{L^{4}(Q)} \le \|\phi\|\|\phi\|_{H^{1}(Q)} \le \|\phi\|^{2}_{H^{1}(Q)} = \|\phi\|^{2} + \|\nabla \phi\|^{2} \le C\|\nabla \phi\|^{2}.
 \end{equation} Using \eqref{I8a}, we infer
 \begin{equation*}
     I_{8} \le c \int_{0}^{t} \|\nabla \Delta \phi_{1}\|^{2}_{L^{4}}\|\bm{u}\|^{2}  + c \int_{0}^{t} \|\nabla \phi\|^{2} .
 \end{equation*} Collecting the above estimates and choosing $\epsilon$ appropriately, we have
 \begin{equation}
 \begin{aligned}
    &\frac{1}{2}\|\bm{u}(t)\|^{2} + c \int_{0}^{t} \|\nabla \bm{u}\|^{2} \le \frac{1}{2}\|\bm{u}(0)\|^{2} + \epsilon \int_{0}^{t} \|\Delta^{2} \phi\|^{2}  \notag \\ \quad& + c \int_{0}^{t} \| \nabla \phi\|^{2} \left( c\ln\left( \frac{c}{\|\nabla \phi\|^{2}} \right)\|D\bm{u}_{2}\|^{2}  +  \|\nabla \Delta \phi_{1}\|^{4}_{L^{4}(Q)}
    + \|\nabla \Delta \phi_{1}\|^{\frac{12}{7}}_{L^{4}(Q)}\right)   \notag \\\quad&
    + c \int_{0}^{t} \|\nabla \phi\|^{2} \left(1 + \|\nabla \phi_{2}\|^{6}_{L^{\infty}(Q)} \right) \\ &+  c\int_{0}^{t} \| \bm{u}\|^{2} \left\{ \|\nabla \bm{u}_{1}\|^{2} + \|\nabla \Delta \phi_{1}\|^{\frac{12}{7}}_{L^{4}(Q)} + \|\nabla \Delta \phi_{1}\|^{2}_{L^{4}(Q)}   \right\} ,
\end{aligned}
 \end{equation} On account of the regularity of $(\bm{u}_{i}, \phi_{i})$, $i=1,2$, the above estimate becomes
 \begin{equation}
\begin{aligned}
    &\frac{1}{2}\|\bm{u}(t)\|^{2} + c \int_{0}^{t} \|\nabla \bm{u}\|^{2}  \le \frac{1}{2}\|\bm{u}(0)\|^{2} + \epsilon \int_{0}^{t} \|\Delta^{2} \phi\|^{2}   \\ \quad& + \int_{0}^{t} R(s)  \| \nabla \phi\|^{2} \left( \ln\left( \frac{c}{\|\nabla \phi\|^{2}} \right) + 1 \right) +  \int_{0}^{t}  R(s)\| \bm{u}\|^{2}, \label{U1}
\end{aligned}
 \end{equation} where $\epsilon > 0$ can be chosen arbitrarily small.

\subsection{The PFC energy identity}
{\color{black} Next, we turn our attention to the second equation in the weak formulation. We have for $i=1,2$ that
\begin{equation} \label{x2}
    ( \phi_{i}', \rho ) + ( \bm{u}_{i} \nabla \phi_{i}, \rho ) + ( m(\phi_{i})\nabla \psi_{i}, \nabla  \rho ) =0, ~~  \forall \hspace{3pt} \rho \in \Phi_1.
\end{equation}
Subtracting \eqref{x2} with $i=2$ from \eqref{x2} with $i=1$, we find that
\begin{equation*}
    \begin{aligned}
        (\phi ' , \rho) + (\bm{u}_{1} \nabla \phi_{1}, \rho)- (\bm{u}_{2} \nabla \phi_{2}, \rho) + (m(\phi_{1})\nabla \psi_{1}, \nabla \rho) - (m(\phi_{2})\nabla \psi_{2}, \nabla \rho) = 0.
    \end{aligned}
\end{equation*}
Thus, taking $\rho = \Delta \phi$ as a test function and integrating in time, we obtain
\begin{align} \label{d2}
    \frac{1}{2}&\|\nabla \phi(t)\|^{2} - \int_{0}^{t} ((\bm{u} \cdot \nabla \phi_{1}), \Delta \phi) - \int_{0}^{t} ((\bm{u}_{2} \cdot \nabla \phi), \Delta \phi ) \notag \\ & + \int_{0}^{t} (m(\phi_{1})\nabla \psi , \nabla \Delta \phi) + \int_{0}^{t} ((m(\phi_{1})-m(\phi_{2}))\nabla \psi_{2}, \nabla \Delta \phi)  = \frac{1}{2}\|\nabla \phi(0)\|^{2},
\end{align}
}where we have set $\psi = \psi_{1} - \psi_{2}$. We can simplify the term involving $\psi$. Indeed, integrating by parts, we get
\begin{align*}
    (m(\phi_{1})\nabla \psi , \nabla \Delta \phi) = - (m'(\phi_{1})(\nabla \cdot \phi_{1}) \psi, \nabla \Delta \phi) - (m(\phi_{1})\psi, \Delta^{2}\phi).
\end{align*}
Observe now that
\begin{align*}
     (m(\phi_{1})\psi, \Delta^{2}\phi) =~& (m(\phi_{1}) (f(\phi_{1}) - f(\phi_{2})), \Delta^{2} \phi) \\ &+ (2m(\phi_{1})\Delta\phi, \Delta^{2} \phi) + (m(\phi_{1}) \Delta^{2} \phi, \Delta^{2}\phi).
\end{align*}
Thus, \eqref{d2} becomes \begin{align*}
    \frac{1}{2}&\|\nabla \phi(t)\|^{2} + \int_{0}^{t} \|\Delta^{2} \phi\|^{2} = \frac{1}{2}\|\nabla \phi(0)\|^{2} +\int_{0}^{t} (\bm{u} \cdot \nabla \phi_{1}, \Delta \phi)  + \int_{0}^{t}  (\bm{u}_{2} \cdot \nabla \phi, \Delta \phi) \\ &+ \int_{0}^{t} (m(\phi_{1})(f(\phi_{1}) - f(\phi_{2})), \Delta^{2} \phi) + 2 \int_{0}^{t} (m(\phi_{1}) \Delta \phi, \Delta^{2} \phi)   \\ &+ \int_{0}^{t} (m(\phi_{1})\Delta^{2}\phi, \Delta^{2}\phi) - \int_{0}^{t} (m'(\phi_{1}) \psi (\nabla \cdot  \phi_{1}) , \nabla \Delta \phi) \\ &- \int_{0}^{t} ((m(\phi_{1})-m(\phi_{2}))\nabla \psi_{2}, \nabla \Delta \phi) ,
\end{align*}
which can be written as
\begin{align}
\label{d3}
     \frac{1}{2}\|\nabla \phi(t)\|^{2} + \int_{0}^{t} \|\Delta^{2} \phi\|^{2}  &= \frac{1}{2}\|\nabla \phi(0)\|^{2} + \sum_{j=9}^{15}I_{j}, ~\text{ for all } t\in [0,T].
\end{align}
We now estimate $I_{j}$ for $j=9,\dots,15$.
\\[1mm]
$\bullet$	 $I_{9}$ $\bullet$
Recalling \eqref{Del} and using the H\"{o}lder and the Young inequalities, we get
\begin{align*}
I_{9} &= \int_{0}^{t} (\bm{u} \cdot \nabla \phi_{1}, \Delta \phi) \mathrm{d}t \le \int_{0}^{t} \|\nabla \phi_{1}\|_{L^{\infty}(Q)} \|\bm{u}\| \|\Delta\phi\|  \\ \quad& \le c \int_{0}^{t} \|\bm{u}\|^{2} \|\nabla \phi_{1}\|^{2}_{L^{\infty}(Q)}  + c \int_{0}^{t} \|\nabla \phi\|^{2}  + \epsilon \int_{0}^{t} \|\Delta^{2}\phi\|^{2}.
 \end{align*}
$\bullet$	 $I_{10}$ $\bullet$
Integrating by parts and using H\"{o}lder's inequality, we have
\begin{align*} I_{10} &= \int_{0}^{t} (\bm{u}_{2} \cdot \nabla \phi, \Delta \phi) = - \int_{0}^{t} (\phi \bm{u}_{2}, \nabla \Delta \phi) \le \int_{0}^{t} \|\bm{u}_{2}\|_{L^{3}(Q)} \|\phi\|_{L^{6}(Q)} \|\nabla \Delta \phi\| \\ \quad& \le c \int_{0}^{t} \|\bm{u}_{2}\|^{2}_{L^{3}(Q)} \|\phi\|^{2}_{H^{1}(Q)} +  \int_{0}^{t} \|\nabla \Delta \phi\|^{2},
 \end{align*} where we also used a Sobolev embedding and Young's inequality. Next, using  \eqref{Nab} and recalling Poincar\'{e}-Wirtinger's inequality, we obtain
 \begin{align*}
     I_{10} &\le c\int_{0}^{t} \|\bm{u}_{2}\|^{2}_{L^{3}(Q)} \|\nabla\phi\|^{2} + \int_{0}^{t} \left( \epsilon \|\Delta^{2}\phi\|^{2} + c\|\nabla \phi \|^{2}\right) \\[2mm] \quad& \le \int_{0}^{t} (c+\|\bm{u}_{2}\|^{2}_{L^{3}(Q)}) \|\nabla \phi\|^{2} + \epsilon \int_{0}^{t} \|\Delta^{2}\phi\|^{2}.
 \end{align*}
$\bullet$	 $I_{11}$ $\bullet$
Recalling that $\phi_{i} \in L^{\infty}(0,T; L^{\infty}(Q))$, $i=1,2$, and using the H\"{o}lder and the Young inequalities, we have
\begin{align*}
 |( (f(\phi_{1}) - f(\phi_{2}), \Delta^{2} \phi)| &\le c(|\phi|, |\Delta^{2}\phi|) \\ \quad  &\le c\|\phi\| \|\Delta^{2}\phi\| \le \epsilon \|\Delta^{2}\phi\|^{2} + c \|\nabla \phi\|^{2}.
 \end{align*}
 Therefore, \begin{equation*}
     I_{11} = \int_{0}^{t} m(\phi_{1})( (f(\phi_{1}) - f(\phi_{2})), \Delta^{2}\phi)  \le \epsilon \int_{0}^{t} \|\Delta^{2}\phi\|^{2}  + c \int_{0}^{t} \|\nabla \phi\|^{2}.
   \end{equation*}
$\bullet$	 $I_{12}$ $\bullet$
Using Young's inequality and the estimate $\|\Delta \phi\|^{2} \le \|\nabla \phi\| \| \nabla \Delta \phi\|$ as well as $\|\phi\|_{H^{4}(Q)} \le c\|\Delta^{2}\phi\|$, it holds
\begin{align*} I_{12} &\le c \int_{0}^{t} \|\Delta \phi\| \|\Delta^{2}\phi\| \le c\int_{0}^{t}  \|\nabla \phi\| \| \nabla \Delta \phi\| + \frac{\epsilon}{2} \int_{0}^{t} \|\Delta^{2}\phi\|^{2} \\ &\le c \int_{0}^{t} \|\nabla \phi\| \|\phi\|_{H^{4}(Q)} + \frac{\epsilon}{2} \int_{0}^{t} \|\Delta^{2}\phi\|^{2} \le c \int_{0}^{t} \|\nabla \phi\|^{2} + \epsilon\int_{0}^{t} \|\Delta^{2}\phi\|^{2}.
 \end{align*}
$\bullet$	 $I_{13}$ $\bullet$
A simple application of \eqref{Nab} and assumption (A1) yields
\begin{align*} I_{13} &\le c \int_{0}^{t} \|\nabla \Delta \phi\|^{2}  \le  c\int_{0}^{t} \|\nabla \phi\|^{2}  + 4\epsilon \int_{0}^{t} \|\Delta^{2}\phi\|^{2}.
 \end{align*}
 $\bullet$ $I_{14}$ $\bullet$ Using assumption (A1), H\"{o}lder's inequality, the estimate $\|\psi\| \le c\|\Delta^{2}\phi\|$ and \eqref{Nab}, one obtains
 \begin{align*}
     I_{14} &= \int_{0}^{t} (m'(\phi_{1}) \psi (\nabla \cdot  \phi_{1}) , \nabla \Delta \phi)  \le c \int_{0}^{t} \|\psi\| \|\nabla \phi_{1} \|_{L^{\infty}(Q)} \|\nabla \Delta \phi\| \\[1ex] &\le c \int_{0}^{t} \|\Delta^{2} \phi\| \|\nabla \Delta \phi \| \le \epsilon \int_{0}^{t} \|\Delta^{2}\phi\|^{2} + c\int_{0}^{t}\|\nabla \phi \|^{2}.
 \end{align*}  $\bullet$ $I_{15}$ $\bullet$
 Using the Lipschitz continuity of $m$ and $\phi_{i}$, we have $m(\phi_{1})-m(\phi_{2}) \le c|\phi_{1}-\phi_{2}|$. Additionally, recalling that $\|\phi\|_{L^{4}(Q)} \le c \| \nabla \phi\|$, $\|\nabla \Delta \phi\| \le c\|\Delta^{2}\phi\|$ and using Young's inequality, we have \begin{align*}
     \int_{0}^{t} ((m(\phi_{1})&-m(\phi_{2}))\nabla \psi_{2}, \nabla \Delta \phi) \le c\int_{0}^{t} \|\phi\|_{L^{4}(Q)} \|\nabla \psi_{2}\| \|\nabla \Delta \phi\|_{L^{4}(Q)} \\[1ex] &\le c\int_{0}^{t} \|\nabla \phi\|\|\nabla \psi_{2}\| \| \Delta^{2} \phi \| \le \epsilon \int_{0}^{t} \|\Delta^{2}\phi\| + c\int_{0}^{t} \|\nabla \phi\|^{2}\|\nabla \psi_{2}\|^{2}.
 \end{align*}
\subsection{Completing the proof of uniqueness}
Collecting the estimates of $I_{9},\dots, I_{15}$, we find
\begin{align*}
    \frac{1}{2}&\|\nabla \phi(t)\|^{2} + \int_{0}^{t} \|\Delta^{2} \phi\|^{2} \le \frac{1}{2}\|\nabla \phi(0)\|^{2} +  c\int_{0}^{t} \|\nabla \phi\|^{2}(1+\|\nabla \psi_{2}\|^{2}+\|\bm{u}_{2}\|^{2}_{L^{3}(Q)}) \\ \quad&+ \epsilon \int_{0}^{t} \|\Delta^{2}\phi\|^{2}   + c \int_{0}^{t} \|\bm{u}\|^{2}\|\nabla \phi_{1}\|^{2}_{L^{\infty}(Q)}.
\end{align*}
Choosing $\epsilon, \epsilon'$ appropriately and using the regularity of $(\bm{u}_{i}, \phi_{i})$ gives us
\begin{equation}
\begin{aligned}
    \frac{1}{2}&\|\nabla \phi(t)\|^{2} + c\int_{0}^{t} \|\Delta^{2}\phi\|^{2}~ \mathrm{d}s \\
    \quad &\le \frac{1}{2}\|\nabla \phi(0)\|^{2} + \int_{0}^{t} R(s)\|\nabla \phi\|^{2}~ \mathrm{d}s +  \int_{0}^{t}R(s)\|\bm{u}\|^{2}~ \mathrm{d}s. \label{U2}
\end{aligned}\end{equation}
Adding \eqref{U1} to \eqref{U2} gives us
\begin{equation}
    \begin{aligned}
    \frac{1}{2}&\|\bm{u}(t)\|^{2} + \frac{1}{2}\|\nabla \phi(t)\|^{2} + c \int_{0}^{t} \|\nabla \bm{u}(s)\|^{2} \mathrm{d}s + c\int_{0}^{t} \|\Delta^{2}\phi\|^{2}\mathrm{d}s \\[1ex] &\le \frac{1}{2}\|\bm{u}(0)\|^{2} + \frac{1}{2}\|\nabla \phi(0)\|^{2} +  \int_{0}^{t} R(s) \| \nabla \phi\|^{2} \left(c \ln\left( \frac{c}{\|\nabla \phi\|^{2}} \right) + \|\nabla \psi_{2}\|^{2} \right) ~ \mathrm{d}s   \\
    &\quad + \int_{0}^{t} R(s) \| \bm{u}(s)\|^{2}~ \mathrm{d}s.    \label{UB}
\end{aligned}
\end{equation}
Observe that \begin{equation*}
    1 \le \ln\left(e\frac{\|\phi\|_{H^{2}}}{\|\nabla \phi\|}\right) \le \ln\left(\frac{\kappa}{\|\nabla \phi\|}\right)
\end{equation*} where $\kappa >0$ is independent of time. Therefore, we can say that \begin{align*}
    \int_{0}^{t}R(s) \| \nabla \phi\|^{2} &\left(c \ln\left( \frac{\kappa}{\|\nabla \phi\|^{2}} \right) + \|\nabla \psi_{2}\|^{2} \right) ~ \mathrm{d}s \\[1ex] &\le c \int_{0}^{t} \|\nabla \phi\|^{2}\ln\left(\frac{\kappa}{\|\nabla \phi\|^{2}}\right) (R(s) + R(s)\|\nabla \psi_{2}\|^{2} ) \mathrm{d}s.
\end{align*}
Recall that $R$ represents a continuous function of time and so $G := R + R\|\nabla \psi_{2}\|^{2}$ in particular belongs to $L^{1}(0,T)$, since $\psi_{2} \in L^{2}(0,T; \Phi_{1})$. Noting that $\|\bm{u}(0)\| =  0 = \|\nabla \phi(0)\|$, we obtain
\begin{align*}
    \frac{1}{2}\|\bm{u}(t)\|^{2} &+ \frac{1}{2}\|\nabla \phi(t)\|^{2} + a_{1} \int_{0}^{t} \|\nabla \bm{u}(s)\|^{2} ~\mathrm{d}s + a_{2} \int_{0}^{t} \| \Delta^{2} \phi(s)\|^{2}~ \mathrm{d}s  \\ \le
    \quad& c\int_{0}^{t} \| \nabla \phi\|^{2}\ln\left(\frac{\kappa}{\|\nabla \phi\|^{2}} \right) G(s)~ \mathrm{d}s + \int_{0}^{t} \|\bm{u}(t)\|^{2} G(s) ~\mathrm{d}s.
\end{align*}
Additionally defining $F(t):= \frac{1}{2}\|\bm{u}(t)\|^{2} + \frac{1}{2}\|\nabla \phi(t)\|^{2}$, the above inequality implies that
\begin{align*}
    F(t) &\le c \int_{0}^{t} \left(\|\nabla \phi\|^{2}\ln\left(\frac{\kappa}{\|\nabla \phi\|^{2}} \right) + \|\bm{u}(s)\|^{2} \right) G(s) ~\mathrm{d}s.
\end{align*} Notice that $x\ln(\frac{\kappa}{x})$ is an increasing function for $x < \frac{\kappa}{e}$, so there exists $C>0$ independent of $t$ such that \begin{equation*}
\|\nabla \phi\|^{2} \ln\left(\frac{\kappa}{\|\nabla \phi\|^{2}}\right) \le F(t)\ln\left(\frac{C}{F(t)}\right),
\end{equation*}
almost everywhere in $[0,T]$. Thus, we get \begin{equation*}
    F(t) \le c \int_{0}^{t} F(s)\ln\left(\frac{C}{F(s)} \right)G(s) ~\mathrm{d}s.
\end{equation*}
Defining $W(s) := s\ln\left(\frac{C}{s}\right)$ for $s \ge 0$, the above inequality can be expressed as
\begin{equation}
\label{ux}
     F(t) \le c \int_{0}^{t} W\left(F(s)\right)G(s) ~\mathrm{d}s.
\end{equation}
Thus, taking $C$ sufficiently large, we can apply Osgood's lemma (see, for instance, \cite[Appendix B]{JFAGiorg}) and deduce that $F(t) = 0$ for almost every $t \in [0,T]$. This implies that $\|\bm{u}(t)\| = \| \nabla \phi(t)\| = 0$ for almost any $t\in [0,T]$, from which we conclude that $(\bm{u}_{1}, \phi_{1}, \psi_{1}) = (\bm{u}_{2}, \phi_{2}, \psi_{2})$ almost everywhere in $Q \times [0,T]$, as required.

Finally, observe that the energy identity \eqref{enineq2} can be recovered through a standard argument, that is, by taking $\bm{v} = \bm{u}$, $\rho = \phi$ as test functions in \eqref{deboleteorema1}, \eqref{deboleteorema2}, respectively (see, for instance,
\cite[Chapt.V, Proof of Prop.V.1.7]{Boyer} for details). The energy identity also implies that, for any $\delta\in (0,T)$ we can find $\tau\in(0,\delta)$ such that
$\bm{u}(\tau)\in\mathbb{V}$ and $\phi(\tau)\in \Phi_3$. Thus, on account of Theorem \ref{theoremstrong},
the solution becomes strong from $\tau$ on.

\begin{remark}
From \eqref{ux} it is also possible to recover a continuous dependence estimate (see \cite[Proof of Theorem 3.3]{JFAGiorg}).
\end{remark}

\section{Proof of Theorem \ref{theoremstrong}}
\label{strong2D}
This section is devoted to the proof of Theorem \ref{theoremstrong}. In Subsection \ref{estimatesstrong}, we obtain stronger a priori estimates on the weak solutions. This improved regularity together with the hypothesis on the initial data allows us to deduce the existence of strong solutions. In Subsection \ref{uniqueness}, we prove that the strong solution depends continuously on the initial data.

\subsection{A priori estimates}
\label{estimatesstrong}
Here, we consider once more the finite dimensional approximation \eqref{eq:Ia}-\eqref{eq:IIIa}.
Exploiting the smoothness of the eigenfunctions, we first establish higher order estimate for the solutions.

$\bullet$	 {\it Step 1:}
Take $\bm{w}=-\Delta\bm{u}_n$ in equation \eqref{eq:Ia} and integrate by parts. By observing that $\nabla \cdot (-\Delta \bm{u}_n) = 0$, we have
\begin{align*}
    &\frac{1}{2}\frac{\d}{\d t} \lVert \nabla\bm{u}_n\rVert_{\mathbb{H}}^2 + \int_Q \eta'(\phi_n)(\nabla \phi_n \cdot \bm{D}\bm{u}_n) \cdot \Delta\bm{u}_n + \int_Q \eta(\phi_n) \vert \Delta \bm{u}_n\vert^2 \nonumber\\
    & -\int_Q (\bm{u}_n \cdot \nabla\bm{u}_n)\cdot \Delta \bm{u}_n + M \int_Q \psi_n \nabla \phi_n \cdot \Delta \bm{u}_n = 0,
\end{align*}
from which, by using Assumption (A1) and the H\"older inequality, we get
\begin{align}
\label{eq:sI}
    &\frac{1}{2}\frac{\d}{\d t} \lVert \nabla\bm{u}_n\rVert_{\mathbb{H}}^2 + \eta_0 \lVert \Delta \bm{u}_n \rVert_{\mathbb{H}}^2 \nonumber \\[2mm]
    &\leqslant \eta_2 \lvert \nabla \phi_n\cdot D\bm{u}_n \rvert_{\mathbb{H}} \lvert \Delta\bm{u}_n \rvert_{\mathbb{H}}  + \lvert \bm{u}_n \cdot \nabla\bm{u}_n \rvert_{\mathbb{H}} \lvert \Delta \bm{u}_n \rvert_{\mathbb{H}} + M \rvert \psi_n \nabla \phi_n \rvert_{\mathbb{H}} \lvert \Delta \bm{u}_n \rvert_{\mathbb{H}}.
\end{align}
We consider each term on the right-hand side of inequality \eqref{eq:sI} separately. By using the Young, the H\"older, and the Ladyzhenskaya inequalities, the first term can be estimated as follows
\begin{align*}
    &\eta_2 \lvert \nabla \phi_n\cdot \bm{D}\bm{u}_n \rvert_{\mathbb{H}} \lvert \Delta\bm{u}_n \rvert_{\mathbb{H}}\leqslant c \lVert \nabla \phi_n \rVert_{\mathbb{H}} \lVert \nabla \phi_n \rVert_{\mathbb{V}} \lVert \bm{D} \bm{u}_n \rVert_{\mathbb{H}} \lVert \bm{D} \bm{u}_n \rVert_{\mathbb{V}}+\frac{\eta_0}{2} \lVert \Delta \bm{u}_n \rVert_{\mathbb{H}}^2,
\end{align*}
from which, by using again the Young inequality and the Korn inequality, we get
\begin{align}
\label{eq:sei}
    &\eta_2 \lvert \nabla \phi_n\cdot D\bm{u}_n \rvert_{\mathbb{H}} \lvert \Delta\bm{u}_n \rvert_{\mathbb{H}} \leqslant c \lVert \phi_n \rVert_{\Phi_2}^4 \lVert \nabla \bm{u}_n \rVert_{\mathbb{H}}^2 + \frac{\eta_0}{8} \lVert \nabla \bm{u}_n \rVert_{\mathbb{H}}^2+\frac{3}{4}\eta_0 \lVert \Delta \bm{u}_n \rVert_{\mathbb{H}}^2.
\end{align}
Similarly, by using the Young, the H\"older, and the Ladyzhenskaya inequalities for the second term, we get
\begin{align*}
    &\lvert \bm{u}_n \cdot \nabla\bm{u}_n \rvert_{\mathbb{H}} \lvert \Delta \bm{u}_n \rvert_{\mathbb{H}}\leqslant  c \lVert \bm{u}_n \rVert_{\mathbb{H}} \lVert \bm{u}_n \rVert_{\mathbb{V}} \lVert \nabla \bm{u}_n \rVert_{\mathbb{H}} \lVert \nabla \bm{u}_n \rVert_{\mathbb{V}} + \frac{\eta_0}{8} \lVert \Delta \bm{u}_n \rVert_{\mathbb{H}}^2,
\end{align*}
which, by means of the Young inequality, entails
\begin{align}
\label{eq:seii}
    &\lvert \bm{u}_n \cdot \nabla\bm{u}_n \rvert_{\mathbb{H}} \lvert \Delta \bm{u}_n \rvert_{\mathbb{H}}\leqslant c \lVert \bm{u}_n \rVert_{\mathbb{H}}^2 \lVert \bm{u}_n \rVert_{\mathbb{V}}^2 \lVert \nabla \bm{u}_n \rVert_{\mathbb{H}}^2  + \frac{\eta_0}{16} \lVert \nabla \bm{u}_n \rVert_{\mathbb{H}}^2 +\frac{3}{16}\eta_0 \lVert \Delta \bm{u}_n \rVert_{\mathbb{H}}^2.
\end{align}
Eventually, by using the Young inequality to the third term, we get
\begin{equation}
\label{eq:seiii}
      M \rvert \psi_n \nabla \phi_n \rvert_{\mathbb{H}} \lvert \Delta \bm{u}_n \rvert_{\mathbb{H}} \leqslant c \lVert \psi_n \nabla \phi_n \rVert_{\mathbb{H}}^2 + \frac{\eta_0}{32} \lVert \Delta \bm{u}_n \rVert_{\mathbb{H}}^2.
\end{equation}
We use \eqref{eq:sei}-\eqref{eq:seiii} in \eqref{eq:sI} and we integrate over $(0,t)$ for $t \in (0,T_n)$, with $T_n < T$, to get
\begin{align}
\label{eq:stI}
   &\frac{1}{2} \lVert \nabla\bm{u}_n(t)\rVert_{\mathbb{H}}^2 + \frac{\eta_0}{32} \int_0^t \lVert \Delta \bm{u}_n \rVert_{\mathbb{H}}^2 \nonumber\\ \leqslant& \;c \int_0^t \lVert \phi_n \rVert_{\Phi_2}^2 \lVert \nabla \bm{u}_n \rVert_{\mathbb{H}}^2 + c \int_0^t \lVert \bm{u}_n \rVert_{\mathbb{H}}^2 \lVert \bm{u}_n \rVert_{\mathbb{V}}^2 \lVert \nabla \bm{u}_n \rVert_{\mathbb{H}}^2+  \frac{3}{16}\eta_0 \int_0^t \lVert \nabla \bm{u}_n \rVert_{\mathbb{H}}^2\nonumber\\ &+ c \int_0^t \lVert \psi_n \nabla \phi_n \rVert_{\mathbb{H}}^2 + \frac{1}{2} \lVert \nabla\bm{u}_n(0)\rVert_{\mathbb{H}}^2.
\end{align}
We now estimate the right-hand side of \eqref{eq:stI}. By first observing that $\lVert \nabla \bm{u}_n(0) \rVert_{\mathbb{H}}^2 \leqslant c$ and by then using the estimates of Subsection \ref{estimatesweak}, we readily get the following bound
\begin{equation}
\label{eq:stII}
\frac{1}{2} \lVert \nabla \bm{u}_n(t) \rVert_{\mathbb{H}}^2 + \frac{\eta_0}{32}\int_0^t \lVert \Delta \bm{u}_n \rVert_{\mathbb{H}}^2 \leqslant c \biggl( 1 + \int_0^t \lVert \bm{u}_n \rVert_{\mathbb{V}}^2 \lVert \nabla \bm{u}_n \rVert_{\mathbb{H}}^2 \biggr).
\end{equation}
As $\bm{u}_n$ is a solution to problem \eqref{eq:Ia}-\eqref{eq:IIIa}, the map $s \mapsto \lvert \bm{u}_n(s) \rvert_{\mathbb{H}^1(Q)}^2 \in C^1([0,T])$ for each $n$. Thus, by applying the Gronwall lemma in \eqref{eq:stII}, we obtain
\begin{equation}
\label{eq:uns1}
\lVert \bm{u}_n \rVert_{L^{\infty}(0,T;\mathbb{V})} \leqslant c,
\end{equation}
which, on account of \eqref{eq:stII}, implies
\begin{equation}
\label{eq:un1}
    \lVert \bm{u}_n \rVert_{L^2(0,T;\mathbb{H}^2(Q))} \leqslant c.
\end{equation}

{\color{black}
At this stage, we bound the remaining terms in equation \eqref{eq:Ia}.
In particular, the inertia term can be handled by using \eqref{eq:uns1} and \eqref{eq:un1} as follows
\begin{equation}
\label{eq:sB0}
    \lVert \bm{u}_n \cdot \nabla \bm{u}_n \rVert_{L^2(0,T;\mathbb{H})}^2 \leqslant \lVert \bm{u}_n \rVert_{L^{\infty}(0,T;\mathbb{H})}^2 \lVert \nabla \bm{u}_n \rVert_{L^2(0,T;\mathbb{L}^{\infty}(Q))}^2 \leqslant c .
\end{equation}
}Moreover, from the estimates \eqref{106}, \eqref{eq:uns1}, \eqref{eq:un1}, and Assumption (A1), we deduce
\begin{align}
\label{a0boundStrong}
    &\lVert\nabla \cdot ( \eta(\phi_n)\bm{D}\bm{u}_n) \rVert_{L^2(0,T;\mathbb{H})} \nonumber\\ &\leqslant c \lVert \eta' \rVert_{L^{\infty}((0,T)\times Q)} \lvert \nabla \phi_n \rvert_{L^{2}(0,T;L^{\infty}(Q))}\lvert \bm{D}\bm{u}_n \rvert_{L^{\infty}(0,T;\mathbb{H})} + c \lvert \Delta \bm{u}_n \rvert_{L^2(0,T;\mathbb{H})}^2 \leqslant c.
\end{align}
Eventually, by comparison in equation \eqref{eq:Ia}, thanks to the estimates \eqref{estimateR}, \eqref{eq:sB0}, and \eqref{a0boundStrong}, we obtain
\begin{equation}
\label{eq:sun'}
    \lVert \bm{u}_n' \rVert_{L^2(0,T;\mathbb{H})} \leqslant c,
\end{equation}
where we also used $\lVert P_{\mathbb{V}^n} \rVert_{\mathcal{L}(\mathbb{V}^n,\mathbb{V}^n)} \leq 1$.

$\bullet$	 {\it Step 2:}
We consider $\rho = -\Delta \psi_n$ in equation \eqref{eq:IIa}. By integrating by parts and by exploiting equation \eqref{eq:IIIa}, we get
\begin{align*}
    &\int_Q \nabla \phi_n'\cdot \nabla (\Delta^2 \phi_n + 2\Delta \phi_n + P_{\Phi_2^n}f(\phi_n)) + \int_Q m(\phi_n)\Delta \psi_n \Delta \psi_n \nonumber\\
    &-\int_Q (\bm{u}_n \cdot \nabla \phi_n) \Delta \psi_n + \int_Q m'(\phi_n)(\nabla \phi_n \cdot \nabla \psi_n) \Delta \psi_n = 0,
\end{align*}
from which, by using Assumption (A1) and the H\"older and the Young inequalities, we obtain
\begin{align*}
    &\int_Q \nabla \phi_n' \cdot \nabla (\Delta^2 \phi_n + 2\Delta \phi_n + P_{\Phi_2^n}f(\phi_n)) +\frac{m_0}{4}\lVert \Delta \psi_n \rVert_{L^2(Q)}^2 \\
    &\leqslant \frac{1}{2m_0} \lVert \bm{u}_n \cdot \nabla \phi_n \rVert_{L^2(Q)}^2 + \frac{m_2^2}{m_0} \lVert \nabla \phi_n \cdot \nabla \psi_n \rVert_{L^2(Q)}^2.
\end{align*}
By integrating by parts the first and the third term on the left-hand side and by observing that $P_{\Phi_2^n} \phi_n' = \phi_n'$ as $ \phi_n' \in \Phi_2^n$, we get
\begin{align}
\label{eq:sIt}
&\frac{\d}{\d t}\biggl( \frac{1}{2} \lVert \nabla \Delta \phi_n \rVert_{L^2(Q)}^2 - \lVert \Delta \phi_n \rVert_{L^2(Q)}^2 \biggr) +\frac{m_0}{4}\lVert \Delta \psi_n \rVert_{L^2(Q)}^2 \nonumber\\
&\leq \frac{1}{2m_0} \lVert \bm{u}_n \cdot \nabla \phi_n \rVert_{L^2(Q)}^2 + \frac{m_2^2}{m_0} \lVert \nabla \phi_n \cdot \nabla \psi_n \rVert_{L^2(Q)}^2 + \int_Q \phi_n' \Delta f(\phi_n).
\end{align}
We consider now each term on the right-hand side of inequality \eqref{eq:sIt} separately. By using the H\"older and the Ladyzhenskaya inequalities, the first term can be estimated as follows
\begin{equation}
\label{eq:se2ii}
    \lVert \bm{u}_n \cdot \nabla \phi_n \rVert_{L^2(Q)}^2 \leqslant c \lVert \bm{u}_n \rVert_{\mathbb{H}} \lVert \bm{u}_n \rVert_{\mathbb{V}} \lVert  \phi_n \rVert_{\Phi_1} \lVert \phi_n \rVert_{\Phi_2}.
\end{equation}
For the second term, by using H\"older, Ladyzhenskaya, and Young inequalities, we readily obtain
\begin{align}
\label{eq:se2iii}
    \frac{m_2^2}{m_0} \lVert \nabla \phi_n \cdot \nabla \psi_n \rVert_{L^2(Q)}^2
    \leqslant& \;c \lVert \nabla \phi_n \rVert_{L^2(Q)}^2 \lVert \nabla\phi_n \rVert_{\Phi_1}^2 \lVert \nabla\psi_n \rVert_{L^2(Q)}^2\nonumber\\[2mm] &+ c \lVert \nabla\psi_n \rVert_{L^2(Q)}^2 + \frac{m_0}{8} \lVert \Delta\psi_n \rVert_{L^2(Q)}^2.
\end{align}
The third term can be treated as follows
\begin{equation}
    \label{eq:se2I}
    \int_Q \phi_n' \Delta f(\phi_n) \leq \frac{1}{2} \lVert \phi_n' \rVert_{\Phi_1'}^2 + \frac{1}{2} \lVert \Delta f(\phi_n) \rVert_{\Phi_1}^2.
\end{equation}
We use \eqref{eq:se2ii}-\eqref{eq:se2I} in \eqref{eq:sIt} and we integrate on $(0,t)$ for $t \in (0,T_n)$, with $T_n < T$, to get
\begin{align}
\label{eq:sItii}
    &\frac{1}{2} \lVert \nabla\Delta\phi_n(t) \rVert_{L^2(Q)}^2  + \frac{m_0}{8} \int_0^t \lVert \Delta \psi_n \rVert_{L^2(Q)}^2  \nonumber \\
    &\leqslant \frac{1}{2} \lVert \nabla \Delta \phi_n(0) \rVert_{L^2(Q)}^2+\lVert \Delta \phi_n(t) \rVert_{L^2(Q)}^2 + \frac{1}{2}\int_0^t\bigl( \lVert \phi_n' \rVert_{\Phi_1'}^2 + \lVert \Delta f(\phi_n) \rVert_{\Phi_1}^2 \bigr)+ c \int_0^t \lVert \psi_n \rVert_{\Phi_1}^2  \nonumber\\
    &\quad + c \int_0^t \lVert \bm{u}_n \rVert_{\mathbb{H}} \lVert \bm{u}_n \rVert_{\mathbb{V}} \lVert \phi_n \rVert_{\Phi_1} \lVert  \phi_n \rVert_{\Phi_2}+ c \int_0^t \lVert \phi_n \rVert_{\Phi_1}^2 \lVert \phi_n \rVert_{\Phi_2}^2 \lVert \psi_n \rVert_{\Phi_1}^2.
\end{align}
We now estimate the right-hand side of \eqref{eq:sItii}. First, we observe that $\lVert \nabla \Delta \phi_n(0) \rVert_{L^2(Q)} \leqslant \lVert \phi_n(0) \rVert_{\Phi_3}$. Thus, thanks to the assumptions on the initial data $\phi_0$ and to the fact that $\lVert P_{\Phi_2^n} \rVert_{\mathcal{L}(\Phi_2^n,\Phi_2^n)} \leq 1$, we have
\begin{equation}
\label{dimstrongphi1}
    \lVert \nabla \Delta \phi_n(0) \rVert_{L^2(Q)}^2  \leqslant c.
\end{equation}
From \eqref{eq:Bphi} we readily obtain
\begin{equation}
\label{tXi}
 \lvert \Delta \phi_n(t) \rvert_{L^2(Q)}^2 \leq \kappa(t), \quad \kappa\in L^{\infty}(0,T).
\end{equation}
Estimates \eqref{phiprimeW} and \eqref{106} imply the uniform bound
\begin{equation}
    \frac{1}{2}\int_0^t \lVert \phi_n' \rVert_{\Phi_1'}^2 + \frac{1}{2}\int_0^t \lVert \Delta f(\phi_n) \rVert_{\Phi_1}^2 \leq c.
\end{equation}
Furthermore, due to the estimates \eqref{eq:Bphi}, \eqref{eq:Bpsi}, and \eqref{eq:uns1} we get
\begin{align}
\label{tempExist}
	 & \int_0^t \lVert \psi_n \rVert_{\Phi_1}^2+ \int_0^t \lVert \bm{u}_n \rVert_{\mathbb{H}} \lVert \bm{u}_n \rVert_{\mathbb{V}} \lVert \phi_n \rVert_{\Phi_1} \lVert  \phi_n \rVert_{\Phi_2}+ \int_0^t \lVert \phi_n \rVert_{\Phi_1}^2 \lVert \phi_n \rVert_{\Phi_2}^2 \lVert \psi_n \rVert_{\Phi_1}^2 \leqslant c.
\end{align}
Putting the estimates \eqref{dimstrongphi1}-\eqref{tempExist} in \eqref{eq:sItii}, we find
\begin{align*}
\frac{1}{2} \lVert \nabla\Delta\phi_n(t) \rVert_{L^2(Q)}^2  + \frac{m_0}{8} \int_0^t \lVert \Delta \psi_n \rVert_{L^2(Q)}^2  \leqslant\; c + \kappa(t),
\end{align*}
from which it is straightforward to deduce
\begin{equation*}
    \lVert \nabla \Delta \phi_n \rVert_{L^{\infty}(0,T;L^2(Q))} + \lVert \Delta \psi_n \rVert_{L^2(0,T;L^2(Q))} \leqslant c,
\end{equation*}
which, combined with the estimates \eqref{eq:Bphi} and \eqref{eq:Bpsi}, implies
\begin{equation}
\label{desiredestimatephiH3}
    \lVert \phi_n \rVert_{L^{\infty}(0,T;\Phi_3)} + \lVert \psi_n \rVert_{L^2(0,T;\Phi_2)} \leqslant c.
\end{equation}
Concerning the equation \eqref{eq:IIa}, we first observe that, by using the estimates \eqref{eq:uns1}, \eqref{desiredestimatephiH3}, and Assumption (A1), we readily obtain the following bound
\begin{equation}
\label{strongboundb1}
    \lVert \bm{u}_n\nabla\phi_n \rVert_{L^{\infty}((0,T)\times Q)} + \lVert \nabla \cdot ( m(\phi_n) \nabla \psi_n ) \rVert_{L^2((0,T)\times Q)} \leqslant c.
\end{equation}
Then, by comparison in equation \eqref{eq:IIa}, using the estimate \eqref{strongboundb1} and $\lVert P_{\Phi_2^n} \rVert_{\mathcal{L}(\Phi_2^n,\Phi_2^n)} \leq 1$, we get
\begin{equation}
\label{fiprimol2S}
    \lVert \phi_n' \rVert_{L^2(0,T;L^2(Q))} \leqslant  c.
\end{equation}

$\bullet$	 {\it Step 3}
As $\psi_n$ is uniformly bounded in $L^2(0,T;\Phi_2)$ (see \eqref{desiredestimatephiH3}), by comparison in equation \eqref{eq:IIIa}, we readily obtain
\begin{equation}
\label{fil6S}
    \lVert \phi_n \rVert_{L^2(0,T;\Phi_6)} \leqslant c.
\end{equation}
Here, we also used the estimate \eqref{boundfW} together with $\lVert P_{\Phi_2^n} \rVert_{\mathcal{L}(\Phi_2^n,\Phi_2^n)} \leq 1$ and \eqref{106}.

Hence, from \eqref{eq:uns1}, \eqref{eq:un1}, \eqref{eq:sun'}, \eqref{desiredestimatephiH3}, \eqref{fiprimol2S}, and \eqref{fil6S}, we extract not relabeled converging subsequences. In particular, in addition to the convergences obtained in the previous section, we have
\begin{align*}
&\bm{u}_n \overset{\ast}{\rightharpoonup} \bm{u}& &\text{in} \qquad L^{\infty}(0,T;\mathbb{V}), \\
&\bm{u}_n \rightharpoonup \bm{u}& &\text{in} \qquad L^2(0,T;\mathbb{H}^2(Q)),\\
&\bm{u}_n' \rightharpoonup \bm{u}'& &\text{in} \qquad L^{2}(0,T;\mathbb{H}),\\
&\phi_n \overset{\ast}{\rightharpoonup} \phi& &\text{in} \qquad L^{\infty}(0,T;\Phi_3),\\
&\phi_n \rightharpoonup \phi& &\text{in} \qquad L^2(0,T;\Phi_6), \\
&\phi_n' \rightharpoonup \phi'& &\text{in} \qquad L^2(0,T;L^2(Q)),\\
&\psi_n \rightharpoonup \psi& &\text{in} \qquad L^2(0,T;\Phi_2).
\end{align*}
Moreover, by applying the Lions-Magenes theorem (see, e.g., \cite[Chapt.II]{Boyer}), we find that $\bm{u} \in C^0([0,T];\mathbb{V})$ and $\phi \in C^0([0,T];H_p^3(Q))$, which conclude the proof of the existence of strong solutions of the problem \eqref{eq:ns}-\eqref{eq:periodic}.

\subsection{Continuous dependence estimate}
\label{uniqueness}
{\color{black} To prove \eqref{dipcont}, let us consider two strong solutions to problem \eqref{eq:ns}-\eqref{eq:periodic}, say, $(\bm{u}_1,\phi_1)$ and $(\bm{u}_2,\phi_2)$ originating from the initial data $(\bm{u}_{1,0},\phi_{1,0})$ and $(\bm{u}_{2,0},\phi_{2,0})$, respectively.
We set $\bm{u}=\bm{u}_1-\bm{u}_2$, $\phi=\phi_1-\phi_2$, and $\psi=\psi_1-\psi_2$ and we deduce from \eqref{deboleteorema1} and \eqref{deboleteorema2} that, for any $\bm{v} \in \mathbb{V}$ and $\rho \in \Phi_1$, the following equations are satisfied
\begin{align}
    &\int_Q \bm{u}' \cdot \bm{v} + \int_Q \eta(\phi_1) \bm{D} \bm{u}_1 \cdot \nabla \bm{v} - \int_Q \eta(\phi_2) \bm{D} \bm{u}_2 \cdot \nabla \bm{v} + \int_Q (\bm{u}_1 \cdot \nabla \bm{u}_1) \cdot \bm{v} \nonumber\\
    &- \int_Q (\bm{u}_2 \cdot \nabla \bm{u}_2)\cdot \bm{v} + M \int_Q \phi_1 \nabla \psi_1 \cdot  \bm{v} - M\int_Q \phi_2 \nabla \psi_2\cdot \bm{v}= 0,\label{eq:421}
\end{align}
and
\begin{align}
    &\int_Q \phi' \rho + \int_Q \bm{u}_1 \cdot \nabla \phi_1 \, \rho - \int_Q \bm{u}_2 \cdot \nabla \phi_2\, \rho + \int_Q m(\phi_1)\nabla \psi_1 \cdot \nabla \rho \nonumber\\&- \int_Q m(\phi_2)\nabla \psi_2\cdot \nabla \rho=0.\label{eq:424}
\end{align}
In particular, due to the regularity of the strong solutions, by choosing $\bm{v} = \bm{u}$ in \eqref{eq:421}, we get
\begin{align*}
    &\int_Q \bm{u}' \cdot \bm{u} + \int_Q \eta(\phi_1)\bm{D}\bm{u}_1 \cdot \nabla \bm{u} - \int_Q \eta(\phi_2)\bm{D}\bm{u}_2 \cdot \nabla\bm{u}+ \int_Q (\bm{u}_1 \cdot \nabla \bm{u}_1) \cdot \bm{u}  \nonumber\\
    &- \int_Q (\bm{u}_2 \cdot \nabla \bm{u}_2) \cdot \bm{u} +M \int_Q \phi_1 \nabla \psi_1 \cdot  \bm{u} - M\int_Q \phi_2 \nabla \psi_2\cdot \bm{u}= 0,
\end{align*}
from which, after minor manipulations and integration by parts, we get
\begin{align}
&\int_Q \bm{u}' \cdot \bm{u} + \int_Q (\eta(\phi_1)-\eta(\phi_2))\bm{D}\bm{u}_1\cdot \nabla \bm{u} + \int_Q \eta(\phi_2) \vert \bm{D} \bm{u} \vert^2+ \int_Q (\bm{u}\cdot \nabla \bm{u}_1) \cdot \bm{u}  \nonumber\\
& -M \int_Q \psi_1 \nabla \phi_1 \cdot  \bm{u} + M\int_Q \psi_2 \nabla \phi_2\cdot \bm{u}= 0. \label{5.136}
\end{align}
}Moreover, the last two terms can be rewritten by making the terms $\psi_1$ and $\psi_2$ explicit according to equation \eqref{deboleteorema3}, by manipulating the equation, and by integration by parts, as follows
\begin{align*}
    &\int_Q \psi_1 \nabla \phi_1 \cdot  \bm{u} - \int_Q \psi_2 \nabla \phi_2\cdot  \bm{u} \nonumber \\&= \int_Q (\Delta^2 \phi_1 + 2\Delta \phi_1 + f(\phi_1)) \nabla \phi_1 \cdot \bm{u} - \int_Q (\Delta^2 \phi_2 + 2\Delta \phi_2 + f(\phi_2) \nabla \phi_2 \cdot \bm{u}\nonumber\\
    &= \int_Q \Delta^2\phi \nabla \phi_1 \cdot \bm{u} + \int_Q \Delta^2\phi_2 \nabla \phi\cdot \bm{u} + 2\int_Q \Delta\phi \nabla\phi_1 \cdot \bm{u} + 2\int_Q \Delta\phi_2 \nabla\phi \cdot \bm{u}.
\end{align*}
Hence, from \eqref{5.136}, we obtain
\begin{align}
\label{eq:uu2}
    &\frac{1}{2}\frac{\d}{\d t} \lVert \bm{u} \rVert_{\mathbb{H}}^2 + \int_Q (\eta(\phi_1)-\eta(\phi_2))\bm{D}\bm{u}_1\cdot \nabla \bm{u} + \int_Q \eta(\phi_2) \vert \bm{D} \bm{u} \vert^2 + \int_Q(\bm{u}\cdot \nabla \bm{u}_1) \cdot \bm{u} \nonumber\\
    &- M\int_Q \Delta^2\phi \nabla \phi_1 \cdot \bm{u} - M\int_Q \Delta^2\phi_2 \nabla \phi\cdot \bm{u} - 2M\int_Q \Delta\phi \nabla\phi_1\cdot \bm{u}\nonumber\\& - 2M\int_Q \Delta\phi_2 \nabla\phi\cdot \bm{u} = 0.
\end{align}

Similarly, relying again on the regularity of the strong solutions, we choose $\rho = \Delta^2 \phi$ in \eqref{eq:424}, namely,
\begin{align*}
    &\int_Q \phi' \Delta^2\phi + \int_Q \bm{u}_1 \cdot \nabla \phi_1 \Delta^2 \phi - \int_Q \bm{u}_2 \cdot \nabla \phi_2 \Delta^2 \phi + \int_Q m(\phi_1)\nabla \psi_1 \cdot \nabla \Delta^2 \phi \nonumber\\&- \int_Q m(\phi_2)\nabla \psi_2\cdot \nabla \Delta^2 \phi=0,
\end{align*}
from which, after minor manipulations and integration by parts, we get
\begin{align}
\label{eq:uf1}
    &\frac{1}{2}\frac{\d}{\d t} \lVert \Delta \phi \rVert_{L^2(Q)}^2 +\int_Q \bm{u} \cdot \nabla\phi_1 \Delta^2 \phi - \int_Q \bm{u}_2 \cdot \nabla\phi \Delta^2 \phi + \int_Q m(\phi_1)\nabla \psi_1 \cdot \nabla \Delta^2 \phi  \nonumber\\ & - \int_Q m(\phi_2)\nabla \psi_2\cdot \nabla \Delta^2 \phi=0.
\end{align}
The last two terms in equation \eqref{eq:uf1} can be rewritten by making the terms $\psi_1$ and $\psi_2$ explicit according to equation \eqref{deboleteorema3} and by manipulating the resulting equation. Thus, from equation \eqref{eq:uf1}, we get
\begin{align}
\label{eq:uf2}
    &\frac{1}{2}\frac{\d}{\d t} \lVert \Delta \phi \rVert_{L^2(Q)}^2 + \int_Q \bm{u} \cdot \nabla \phi_1 \Delta^2 \phi - \int_Q \bm{u}_2 \cdot \nabla \phi \Delta^2 \phi  \nonumber\\
    &+ \int_Q (m(\phi_1)-m(\phi_2))\nabla \Delta^2 \phi_1 \cdot \nabla \Delta^2 \phi +\int_Q m(\phi_2) \vert\nabla \Delta^2 \phi \vert^2 + 2\int_Q m(\phi_2) \nabla \Delta \phi \cdot \nabla \Delta^2 \phi  \nonumber \\
    &+ 2\int_Q (m(\phi_1)-m(\phi_2))\nabla \Delta \phi_1 \cdot \nabla \Delta^2 \phi + \int_Q m(\phi_2)f'(\phi_2) \nabla \phi \cdot \nabla \Delta^2 \phi \nonumber\\& +\int_Q (m(\phi_1)f'(\phi_1)-m(\phi_2)f'(\phi_2)) \nabla \phi_1 \cdot \nabla \Delta^2 \phi  = 0.
\end{align}
We sum equation \eqref{eq:uu2} divided by $M$ with equation \eqref{eq:uf2} and we obtain
\begin{align*}
    &\frac{1}{2}\frac{\d}{\d t} \lVert \Delta \phi \rVert_{L^2(Q)}^2 + \frac{1}{2M}\frac{\d}{\d t} \lVert \bm{u} \rVert_{\mathbb{H}}^2 + \frac{1}{M} \int_Q (\eta(\phi_1)-\eta(\phi_2))D\bm{u}_1\cdot \nabla \bm{u}  \nonumber\\
    &+ \frac{1}{M} \int_Q \eta(\phi_2) \vert D \bm{u} \vert^2 +\frac{1}{M} \int_Q B_0(\bm{u},\bm{u}_1)\cdot \bm{u} - \int_Q \Delta^2\phi_2 \nabla \phi \cdot \bm{u} - 2\int_Q \Delta\phi \nabla\phi_1 \cdot \bm{u} \nonumber \\
    & - 2\int_Q \Delta\phi_2 \nabla\phi\cdot \bm{u} - \int_Q \bm{u}_2 \cdot \nabla \phi \Delta^2 \phi + \int_Q (m(\phi_1)-m(\phi_2))\nabla \Delta^2 \phi_1 \cdot \nabla \Delta^2 \phi   \nonumber \\
    &+ \int_Q m(\phi_2) \vert\nabla \Delta^2 \phi \vert^2 +2\int_Q (m(\phi_1)-m(\phi_2))\nabla \Delta \phi_1 \cdot \nabla \Delta^2 \phi +2\int_Q m(\phi_2) \nabla \Delta \phi \cdot \nabla \Delta^2 \phi  \nonumber \\
    &+\int_Q (m(\phi_1)f'(\phi_1)-m(\phi_2)f'(\phi_2)) \nabla \phi_1 \cdot \nabla \Delta^2 \phi + \int_Q m(\phi_2)f'(\phi_2) \nabla \phi \cdot \nabla \Delta^2 \phi = 0.
\end{align*}
{\color{black} Thus, by using Assumption (A1) and the Korn inequality, we obtain
\begin{align}
\label{eq:uuf4}
    &\frac{1}{2}\frac{\d}{\d t} \biggl( \lVert \Delta \phi \rVert_{L^2(Q)}^2 + \frac{1}{M} \lVert \bm{u} \rVert_{\mathbb{H}}^2 \biggr) + \frac{\eta_0}{M} \lVert \bm{u} \rVert_{\mathbb{V}}^2 + m_0\lVert \nabla\Delta^2\phi \rVert_{L^2(Q)}^2  \nonumber\\
    &\leqslant \frac{1}{M} \int_Q \vert (\eta(\phi_1)-\eta(\phi_2))D\bm{u}_1\cdot \nabla \bm{u} \vert +\frac{1}{M} \int_Q \vert (\bm{u} \cdot \nabla \bm{u}_1)\cdot \bm{u} \vert + \int_Q  (\nabla \phi \cdot \bm{u})\Delta^2\phi_2  \nonumber\\
    &\quad+ \int_Q  (\nabla \phi  \cdot \bm{u}_2 ) \Delta^2 \phi + 2 \int_Q (\nabla\phi_1 \cdot \bm{u}) \Delta\phi + 2\int_Q  (\nabla\phi \cdot \bm{u})\Delta\phi_2   \nonumber\\
    &\quad+ \int_Q \vert (m(\phi_1)-m(\phi_2))\nabla \Delta^2 \phi_1 \cdot \nabla \Delta^2 \phi \vert +2\int_Q \vert (m(\phi_1)-m(\phi_2))\nabla \Delta \phi_1 \cdot \nabla \Delta^2 \phi \vert   \nonumber\\
    &\quad + 2\int_Q m(\phi_2) \nabla \Delta \phi \cdot \nabla \Delta^2 \phi +\int_Q \vert (m(\phi_1)f'(\phi_1)-m(\phi_2)f'(\phi_2)) \nabla \phi_1 \cdot \nabla \Delta^2 \phi \vert\nonumber\\
    &\quad + \int_Q \vert m(\phi_2)f'(\phi_2) \nabla \phi \cdot \nabla \Delta^2 \phi\vert \nonumber\\
    &=:\sum_{j=1}^{11} I_j.
\end{align}
}We now estimate each term $I_j$ on the right-hand side of \eqref{eq:uuf4} separately. First, we note that, as the solutions $(\bm{u}_1,\phi_1),(\bm{u}_2,\phi_2)\in C^0([0,T];\mathbb{V})\times C^0([0,T];H^3(Q))$, there exists $R \in C^0([0,T])$ such that
\begin{align}
&\lVert \bm{u}_j(t) \rVert_{\mathbb{V}} \leqslant R(t), \label{eq:limitu}\\
&\lVert \phi_j(t) \rVert_{H^3(Q)}\leqslant R(t) , \label{eq:limitf}
\end{align}
for all $t\in [0,T]$ and $j\in\{1,2\}.$

A {\it caveat} on notation: in the following we use the same symbol $R(\cdot)$ to indicate a time-dependent function, possibly depending on $m_0, M_0, m_2, \eta_0, \eta_1, \eta_2 , M$, which belongs to $C^0([0,T])$. Note that $R(\cdot)$ may change, even within the same line. Similarly, we use the same symbol $c$ to indicate a positive constant, possibly depending on $m_0, M_0, m_2, \eta_0, \eta_1, \eta_2 , M$, which may change, even within the same line.

We aim to exploit the Gronwall lemma to establish a uniform bound for the solution. To that end, we focus on the right-hand side of inequality (\ref{eq:uuf4}). Classical inequalities like H\"older, Minkowsky, Young, and Poincar\'e, along with Assumption (A1), the Korn equality, Lemma \ref{poincphi}, and Sobolev Embedding Theorem are employed for the manipulation of the terms $I_j$. This manipulation serves a twofold purpose. Firstly, it yields terms proportional to $\lVert \Delta \phi \rVert_{L^2(Q)}^2$ and $\lVert \bm{u} \rVert_{\mathbb{H}}^2$, which are instrumental in applying the Gronwall lemma. Secondly, it produces terms proportional to $\lVert \bm{u} \rVert_{\mathbb{V}}^2$ and $\lVert \nabla \Delta^2 \phi \rVert_{L^2(Q)}^2$. Notably, the coefficients multiplying these latter terms are strategically chosen to ensure their absorption into the corresponding terms on the left-hand side.

{\color{black} These computations and manipulations yield the following estimates
\begin{align}
    I_1 &\leqslant R(t) \lVert \Delta \phi \rVert_{L^2(Q)}^2 + \frac{\eta_0}{2M} \lVert \bm{u} \rVert_{\mathbb{V}}^2, \label{eq:I1}\\
    I_2 &\leqslant  R(t) \lVert \bm{u} \rVert_{\mathbb{H}}^2 + \frac{\eta_0}{4M} \lVert \bm{u} \rVert_{\mathbb{V}}^2,\\
    I_3 &\leqslant R(t) \lVert\Delta\phi\rVert_{L^2(Q)}^2 +  \frac{\eta_0}{8M}\lVert\bm{u}\rVert_{\mathbb{V}}^2,\\
    I_4 &\leqslant R(t)\lVert\Delta\phi\rVert_{L^2(Q)}^2 +  \frac{m_0}{2}\lVert\nabla\Delta^2\phi\rVert_{L^2(Q)}^2,\\
    I_5 &\leqslant R(t) \bigl( \lVert\bm{u}\rVert_{\mathbb{H}}^2 + \lVert\Delta \phi\rVert_{L^2(Q)}^2 \bigr),\\
    I_6 &\leqslant R(t)\lVert\Delta\phi\rVert_{L^2(Q)}^2 + \lVert\bm{u}\rVert_{\mathbb{H}}^2,\\
    I_7 &\leqslant c \lVert\nabla\Delta^2\phi_1\rVert_{L^2(Q)}^2 +\frac{m_0}{4}\lVert\nabla\Delta^2\phi\rVert_{L^2(Q)}^2,\\
    I_8 &\leqslant R(t) \Vert \Delta \phi \rVert_{L^2(Q)}^2 +\frac{m_0}{8}\lVert\nabla\Delta^2\phi\rVert_{L^2(Q)}^2,\\
    I_9 &\leqslant c \lVert \Delta \phi \rVert_{L^2(Q)}^2 + \frac{m_0}{16}\lVert\nabla\Delta^2\phi\rVert_{L^2(Q)}^2 ,\\
    I_{10} &\leqslant R(t)\lVert \Delta \phi \rVert_{L^2(Q)}^2 + \frac{m_0}{32} \lVert\nabla\Delta^2\phi\rVert_{L^2(Q)}^2,\\
    I_{11} &\leqslant R(t) \lVert\Delta\phi\rVert_{L^2(Q)}^2 + \frac{m_0}{64}\lVert\nabla\Delta^2\phi\rVert_{L^2(Q)}^2.\label{eq:I11}
\end{align}

Taking estimates \eqref{eq:I1}-\eqref{eq:I11} into account, we deduce from \eqref{eq:uuf4} that
\begin{align}
\label{ved}
    &\frac{1}{2}\frac{\d}{\d t} \biggl( \lVert \Delta \phi \rVert_{L^2(Q)}^2 + \frac{1}{M} \lVert \bm{u} \rVert_{\mathbb{H}}^2 \biggr) + \frac{\eta_0}{8 M} \lVert \bm{u} \rVert_{\mathbb{V}}^2 + \frac{m_0}{64}\lVert \nabla\Delta^2\phi \rVert_{L^2(Q)}^2  \nonumber\\[2mm]
    &\leqslant C(t) \biggl( \lvert \Delta \phi \rvert_{L^2(Q)}^2 + \lvert \bm{u} \rvert_{\mathbb{H}}^2 \biggr),
\end{align}
where $C(t) = R(t)\bigr(1 + \lvert  \phi_1 \rvert_{\Phi_5}^2 \bigl)$ is summable in $(0,T)$.
}In particular, estimate \eqref{ved} implies that
\begin{equation*}
    \frac{1}{2}\frac{\d}{\d t} \biggl( \lVert \Delta \phi \rVert_{L^2(Q)}^2 + \frac{1}{M} \lVert \bm{u} \rVert_{\mathbb{H}}^2 \biggr) \leqslant C(t)\biggl( \lVert\Delta\phi\rVert_{L^2(Q)}^2 +\lVert\bm{u}\rVert_{\mathbb{H}}^2\biggr).
\end{equation*}
Observe now that $\lVert\Delta\phi\rVert_{L^2(Q)}^2 +\lVert\bm{u}\rVert_{\mathbb{H}}^2 \in L^{\infty}(0,T)$. Then, integrating over $(0,t)$ with $t\leqslant T$ and applying the Gronwall lemma yield
\begin{equation}
\label{tempoGronw}
    \lVert\Delta\phi(t)\rVert_{L^2(Q)}^2 +\lVert\bm{u}(t)\rVert_{\mathbb{H}}^2 \leqslant c \bigl( \lVert\Delta\phi(0)\rVert_{L^2(Q)}^2 +\lVert\bm{u}(0)\rVert_{\mathbb{H}}^2 \bigr),
\end{equation}
from which, thanks to the regularity of $\phi$ and $\bm{u}$, we conclude
\begin{equation}
\label{dopogronwall}
    \lVert \Delta\phi \rVert_{C^0([0,T];L^2(Q))}^2 + \lVert\bm{u}\rVert_{C^0([0,T];\mathbb{H})}^2 \leqslant c \bigl( \lVert \Delta\phi(0) \rVert_{L^2(Q)}^2 + \lVert\bm{u}(0)\rVert_{\mathbb{H}}^2 \bigr).
\end{equation}
Furthermore, we use \eqref{tempoGronw} in \eqref{ved} to obtain
\begin{align*}
&\frac{1}{2}\frac{\d}{\d t} \biggl( \lVert \Delta \phi \rVert_{L^2(Q)}^2 + \frac{1}{M} \lVert \bm{u} \rVert_{\mathbb{H}}^2 \biggr) + \frac{\eta_0}{8 M} \lVert \bm{u} \rVert_{\mathbb{V}}^2 + \frac{m_0}{64}\lVert \nabla\Delta^2\phi \rVert_{L^2(Q)}^2  \nonumber\\[2mm]
&\leqslant C(t) \bigl( \lVert\Delta\phi(0)\rVert_{L^2(Q)}^2 +\lVert\bm{u}(0)\rVert_{\mathbb{H}}^2 \bigr),
\end{align*}
thus, by integrating over $(0,T)$, we finally get
\begin{equation}
\label{tempGronw2}
    \lVert\bm{u}\rVert_{L^2(0,T;\mathbb{V})}^2 + \lVert\nabla\Delta^2\phi\rVert_{L^2(0,T;L^2(Q))}^2 \leqslant c \biggl( \lVert \Delta\phi(0) \rVert_{L^2(Q)}^2 + \lVert\bm{u}(0)\rVert_{\mathbb{H}}^2 \biggr).
\end{equation}
Eventually, adding together \eqref{dopogronwall} and \eqref{tempGronw2}, we obtain
\begin{align*}
    &\lVert \Delta\phi \rVert_{C^0([0,T];L^2(Q))}^2 + \lVert\bm{u}\rVert_{C^0([0,T];\mathbb{H})}^2 + \lVert\bm{u}\rVert_{L^2(0,T;\mathbb{V})}^2 + \lVert\nabla\Delta^2\phi\rVert_{L^2(0,T;L^2(Q))}^2 \\
    &\leqslant c \biggl( \lVert \phi_{0,1}-\phi_{0,2} \rVert_{\Phi_2}^2 + \lVert\bm{u}_{0,1}-\bm{u}_{0,2}\rVert_{\mathbb{H}}^2 \biggr),
\end{align*}
which implies \eqref{dipcont}. The proof is finished.


\bigskip
{\bf Acknowledgments.} This research was founded in part by the Austrian Science Fund (FWF) grants \href{https://doi.org/10.55776/P32788}{10.55776/P32788} and \href{https://doi.org/10.55776/I5149}{10.55776/I5149}. For open access purposes, the authors have applied a CC BY public copyright license to any author accepted manuscript version arising from this submission.
The authors thank Andrea Giorgini for having made the connection possible. C. Cavaterra and M. Grasselli
are members of Gruppo Nazionale per l'Analisi Matematica, la Probabilit\`{a} e le loro Applicazioni
(GNAMPA), Istituto Nazionale di Alta Matematica (INdAM). C. Cavaterra and M. Grasselli have been partially supported by MIUR-PRIN
Grant 2020F3NCPX ``Mathematics for industry 4.0 (Math4I4)''. C. Cavaterra's work has also been partially supported by MIUR-PRIN Grant 2022 ``Partial differential equations and related geometric-functional inequalities''. Moreover, the research of C. Cavaterra and M. Grasselli
is part of the activities of ``Dipartimento di Eccellenza 2023-2027'' of Universit\`a degli Studi di Milano (C. Cavaterra) and Politecnico di Milano (M. Grasselli).

\end{document}